\documentclass[a4paper,10pt]{article}
\usepackage{hyperref}
\usepackage{amsfonts}
\usepackage{amssymb}
\usepackage[utf8]{inputenc}

\usepackage[english ]{babel}

\usepackage{bbding}

\usepackage[all]{xy}

\usepackage{tikz-cd} 
\usepackage[all]{xy}
\usepackage{amsfonts}
\usepackage{amssymb}
\usepackage{graphicx}
\usepackage{mathtools}
\usepackage{skak} 
\usepackage{foekfont}

\usepackage{calligra}
\usepackage[T1]{fontenc}

\usepackage{amsmath,mathdots}

\usepackage{stmaryrd}

\usepackage[utf8]{inputenc}

\usepackage{mathtools}
\usepackage{mathdots}
\usepackage{tikz}

\definecolor{amber(sae/ece)}{rgb}{1.0, 0.49, 0.0}
\newfont{\rsfsten}{rsfs10 scaled 1200}
\usepackage{MnSymbol}
\usepackage{tikz}
\usepackage{amscd}
\usepackage{MnSymbol}
\usepackage{wasysym}

\usepackage{mathrsfs}

\usepackage{pdfcolmk}

\usepackage{graphicx}

\newcommand*{\rom}[1]{\expandafter\@slowromancap\romannumeral #1@}

\usepackage{wrapfig}

\newcommand{\tightunderset}[2]{%
  \mathop{#2}\limits_{\vbox to .3ex{\kern-0.95ex\hbox{$#1$}\vss}}}
\newcommand{\tightoverset}[2]
{%
  \mathop{#2}\limits_{\vbox to .3ex{\kern-0.95ex\hbox{$#1$}\vss}}}

\newcommand{\oset}[2]{%
  {\mathop{#2}\limits^{\vbox to -.5\ex@{\kern-\tw@\ex@
   \hbox{\scriptsize #1}\vss}}}}
\usepackage[utf8]{inputenc}

\usepackage {ifsym}

\usepackage {url}
\usepackage{rotating}

\usepackage{mathrsfs}
\usepackage{amsbsy}

\title { Curvature, Expansion, Kolmogorov's Diameter, Hilbert's Rational Designs and    Overtwisted Immersions I}

\author{Misha Gromov} 

\begin{document}
\maketitle

\begin
{abstract}

We   prove the existence of  {\it  locally distance increasing maps}  with   {\it controllably small curvatures} 
between Riemannian manifolds,
where our main construction
depends  on  the presence of  {\it particular spherical  and almost
spherical 
sections} of the unit balls in the $l_{p=4}$ spaces.

   In the part II [Gr 2022'] we prove similar results for    $C^\infty$-smooth {\it isometric immersions} $X^m\to Y^N$,
    where our approach   allows an   improvement of
    the present-day bounds on  the  dimension $N$ of the ambient manifold $Y$ in certain cases.

  \end{abstract}
  \tableofcontents


\section {Introduction}

\textbf{ Immersions}  are $C^1$-maps $f:X\to Y$ between smooth manifolds, such that
 their differentials $df:T(X)\to T(Y)$ nowhere vanish,\footnote{Immersions are locally one-to-one but  globally  they may have self intersections. Immersions without self intersections are called 
{\it embeddings}, where, if $X$ is non-compact, one may  require   
  the induced topology in $X$ to be equal the original one.}
$$df(\tau)=0\implies \tau=0, \tau\in T(X).$$
     
    {\it \textbf  \it \textbf {The {\rm (maximal normal bundle)} curvature}} of an immersed    $X$ in a Riemannian  $Y$,
    $$f: X\hookrightarrow  Y, $$
       is
  {\it the supremum of the $Y$-curvatures }  of geodesics $\gamma\subset X$, 
    for the induced Riemannian metric in $X$, 
  $$curv^\perp(X)=curv^\perp(f(X))=curv_f^\perp(X)=curv^\perp(X\overset {f }\hookrightarrow  Y )  =curv^\perp(X\hookrightarrow  Y ),$$

  {  \textbf {Minimal Curvature Problem.}} What is  the 
  infimum of curvatures  of immersions $f:X\hookrightarrow  Y$,
$$min.curv^\perp(X, Y)=min.curv^\perp(X\hookrightarrow Y)?$$

{\it \color {red!33!black} Product Example.}  If $X$ is a product of spheres, 
$$X=\bigtimes _{i=1}^lS^{m_i},$$
 and $Y$  is the unit ball $B^N(1)\subset \mathbb R^N$ then (apart from the trivial case of $l=1$) we  know the exact value of 
  $min.curv^\perp(\bigtimes _{i=1}^lS^{m_i}, B^N(1)$ 
 {\color {magenta}only} where all $m_i=1$, i.e. for the torus $\mathbb T^l$,  and  where
 $N$ is large:
$$min.curv^\perp(\mathbb T^l, B^N(1))=\sqrt {3{l\over l+2}}, \mbox { }   N>>l^2. \leqno {\color {blue} \left[\sqrt 3\right]_{\mathbb T}}$$
  (See  sections 3, 5  and  [Pet2023].) 
 
 But if all  $m_i=2 $, for instance, i.e. $X=(S^2)^l$ we  {\color {red!45!black} neither can show} that   
 $$ min.curv^\perp((S^2)^l,B^{2l+1})\to \infty \mbox  { for }l\to\infty $$ 
 {\color {red!45!black} nor that } 
 $$ \frac {min.curv^\perp((S^2)^l,B^{10l})}{\sqrt l}\to 0 \mbox  { for }l\to\infty. $$ 



\vspace {1mm}



\vspace{3mm}

\hspace{33mm}{\sc Hand-Made Immersions}

\textbf {   Clifford  Embeddings.} The product $X$ of spheres $S^{m_i}(r_i)\subset \mathbb R^{m_i+1}$, 
  $i=1,...,l$, for $\sum_{i=1}^lr_i^2=1$ naturally isometrically  imbeds   to the boundary of the unit $N$-ball for 
  $N=k+\sum_im_i$:
  
  $${\sf Cl} :X=S^{m_1}(r_1)\times ... \times S^{m_l}(r_l)\to S^{N-1}(1)\subset B^N(1)\subset
   \mathbb R^{m_i+1}\times ..\times \mathbb R^{m_i+1}$$
where, clearly, 
$$curv^\perp(X\overset {\sf Cl}\subset B^N)=\max_i1/r_i.$$ 

This, for $r_1=r_2=...=r_l$,  delivers a  codimension $l$-embedding  with  curvature $\sqrt l$. Thus, 
  $$min.curv^\perp\left(\bigtimes_{i=1}^l S^{m_i}, B^N(1)\right )\leq \sqrt l, \mbox { } 
 N=l+\sum_im_i. $$

If $l=1$, then this  is optimal.  In fact, it is  obvious  that    
$$curv\left (X\hookrightarrow B^m(1)\times \mathbb R^N\right)\geq 1,
\mbox{  for  } n\geq   2.$$ 
for all smoothly immersed closed  $m$-manifolds 
 $X$ in  the "unit band" $B^m(1)\times \mathbb R^N$.

But, for instance,  the {\it equality} 
$$min.curv^\perp(\mathbb T^m\hookrightarrow B^{2m})= \sqrt m$$ 
is { \color {red!67!black}problematic   for all $m\geq 2$}.

\textbf { Round  $\bf m$-Tori  in the  Unit  $\bf (m+1)$-Balls.}
$$min.curv^\perp(\mathbb T^2\hookrightarrow B^3)\leq 3:$$
{\sl the boundary of the   $ \frac {1}{3}$-neighbourhood of the circle of radius $\frac {2}{3}$
 in the space  
has  $curv^\perp(\mathbb T^2\subset \mathbb R^3)=3.$} 

 Similarly (see  section 4.1) 
$$min.curv^\perp(\mathbb T^3\hookrightarrow B^4)\leq 2\sqrt2+1<4$$
$$min.curv^\perp(\mathbb T^7=\mathbb T^3\times \mathbb T^3\times \mathbb T^1\hookrightarrow B^8)\leq 8+2\sqrt 2+1<12$$
\hspace {24mm} ...................................................................

$$min.curv^\perp(\mathbb T^m, B^{m+1}) <m^\frac{3}{2}, \mbox { } m=2^k-1.$$

    {\color {teal}\textbf{ Veronese embeddings}}\footnote {These are flashes from a superior world.} of the real projective spaces satisfy (see 5.1),
 $$curv\left (\mathbb RP^m \hookrightarrow  B^\frac {m(m+3)}{2}\right)=\sqrt \frac {2m}{m+1}, \mbox{  e.g.} $$
   $$curv\left (\mathbb RP^2 \hookrightarrow  B^5\right)=2\sqrt \frac {1}{3}<1.155.$$

{\bf Conjecture.} $$min.cirv(X^m,B^N)< \sqrt \frac {2m}{m+1}\implies X=_{diffeo} S^m.$$

\subsection { Immersions with Small Curvature and $\mathcal D(m,N)$-Approximation}\label {examples}


 \textbf { Expansion.} A map between metric spaces, 
$$f:X\to Y,$$ 
 is {\it $\lambda$-expanding, $\lambda>0$}, if it increases  the  the length of  
 curves $\xi:[0:1]\to X$ by a factor $\geq \lambda$, 
$$ length(f\circ \xi)\geq \lambda\cdot length(\xi) \mbox { for all continuous maps $\xi:[0:1]\to X$}. $$
continuous maps. 

{\it Expanding} is an abbreviation for "1-expanding".
 
{\it Riemannian Example.} A $C^1$-{\it smooth} map $f$ between   Riemannian manifolds, e.g.   open subsets in Euclidean spaces,
is $\lambda$-expanding if and only if $||df(\tau)||\geq|| \lambda\tau ||$ for all tangent vectors $\tau\in T(X).$

 Thus, smooth  expanding maps are immersion and every immersion $f$  expands with  respect to  some Riemannian metrics $g=g(f)$ in $X$ and $h=h(f)$ in $Y$.

{\it Equidimensional example.} If $dim(X)=dim(Y)$  then smooth immersions $X\hookrightarrow Y$ are local diffeomorphisms and $C^1$-smooth expanding maps are locally distance increasing. \footnote  {Expanding locally 
homeomorphic maps are also locally distance increasing, but the  absolute value   map $x\mapsto |x|$, for example,  is 1-expanding but not locally homeomorphic.}
  
    {\it \textbf {The relative}} (maximal) {\it \textbf {curvature}} of an immersion between Riemannian  
  manifolds, 
  $$(X,g) \hookrightarrow (Y,h) $$  
  is
   {\it the supremum of $h$-curvatures   in $Y$},  of $g$-geodesics $\gamma\subset X$, 
   $$ curv(f)=curv^X(f)=curv_Y^X(f)=curv_h^g(f)=\sup_{\gamma\subset X}curv_h(f(\gamma)).$$
   
   If $g=f^\ast(h)$ is the induced Riemannian metric in $X$, this is  our   curvature of $X$ in $Y$,
  $$curv_h^g(f)=curv^\perp(X\overset {f }\hookrightarrow  Y ). $$
   (This $curv^\perp (X)$ unlike $curv(f)$  is defined for immersions of smooth manifolds with no metrics on them.)

 {\it Equidimensional example.} 
If $dim(X)=dim(Y)$,  then  $curv^\perp(X\overset {f}\hookrightarrow  Y )=0$, while $curv^X(f)$ measures by how much $f$ deviates from a projective map.

 \textbf {Normal Immersions, where $\mathbf {curv_F^\perp(X)=curv^X(f)}$.} Call an  immersion between Riemannian 
 manifolds $f:X(g) \hookrightarrow Y(h)$  {\it normal} if  for all normal vectors to $X$ in $Y$,
 $$\nu\in T^\perp_x(X)=T_f(x)(Y)\ominus df(T_x(X))$$
 the second  quadratic form II$_\nu$ of the  immersed $X \overset{f}\hookrightarrow$
  is  {\it simultaneously diagonalizable} with the quadratic forms  $g(x)$ and  $f^\ast(h)$ on the tangent space $T_x(X)$.
   For instance,  isometric immersions are normal.

   Clearly, $curv^\perp_f(X)=curv^X(f)$  for {\it  isometric} immersions $f$

\textbf {Curvature in Spheres.} If an immersion $X\to S^{N-1}(1)$ is normal then so is the
 corresponding immersion to  $\mathbb R^N\supset  S^{N-1}(1)$, where
the spherical 
  curvature of $X$   is related to the Euclidean one  by the Pythagorean
 theorem: 
 $$(curv^\perp(X\hookrightarrow  S^{N-1}(1))^2=(curv^\perp(X\hookrightarrow \mathbb R^N)^2-1.$$

Notice that the Clifford embeddings to the unit sphere are known to be {\it optimal}  for  $l=2$,
$$min.curv^\perp(S^{m_1}\times S^{m_2},S^{m_1+m_2+1}(1))=1, \mbox{ } m_1,m_2\geq1,\footnote {See [Ge2021], section  3.7.3 in  [Gr2022] 
 and section 5.5 in the present paper.}$$
 but the  corresponding Euclidean  equality 
 $$min.curv^\perp(S^{m_1}\times S^{m_2},B^{m_1+m_2+2}(1))=\sqrt 2,$$
 remains {\color {red!66!black} conjectural}  for all  
  $m_1,m_2\geq 1$, except   for $m_1=m_2=1$ [Pet].
 
.

\textbf { Curvature in Codimension 1.} This curvature of $X^m\hookrightarrow Y^{m+1}$
 is the supremum of the principal curvatures of $X$ in $Y$ over all points $x\in X$.

Here normality means that the induced quadratic form $f^\ast(g)(x)$ on 
the tangent space $T_x(X)$  is,  at all $\in X$, diagonalizabel in the same basis as the second fundamental form II of $X$.

{\it Example.}  the   immersion   $ S^m(r)\times S^1\to \mathbb R^{m+2}$  obtained by rotating 
  $S^m(r)\hookrightarrow \mathbb R^{m+1}$ around  a line in  $\mathbb R^{m+1}$  within distance $R>r$ from the origin 
is normal  with curvature $\max \left (\frac {1}{R} , \frac {1}{R-r}\right)$.

\vspace {2mm}

 \hspace {11mm} {\sc Expanding Immersions and Regular Homotopies}\vspace {1mm}

The minimal curvature problem  can be refined in two  ways as follows.\vspace {1mm}

{\sf What is the  minimal curvature of {\it expanding  immersions} between given
{\it Riemannian}   manifolds?}

What is the  minimal curvature in a given  homotopy  or regular homotopy 
\footnote {{A $C^1$-continuous  homotopy $f_t$ of smooth maps 
 is {\it regular} if the maps  $f_t$ are {\it immersions} for all $t$.}} class of immersions? \vspace {1mm}

 Below are  partial answers to these questions.\vspace {1mm}

 $\pmb{\mathcal D}(\mathbf {m,N)}$: {\it  \textbf {Curvature of Euclidean Expanding Maps}.   Let  ${\mathcal D}( {m,N)}$} be the infimum of the relative curvatures of the smooth expanding maps $f$ from
  the Euclidean $m$-space to the unit $N$-ball,
  $$ \mathcal D(m,N)=\inf_fcurv^{{\mathbf e}_m}_{{\mathbf e}_N}(f),$$
where $\mathbf e_m$  and $\mathbf e_N$ denote the Euclidean metrics in $\mathbb R^m$ and
 $\mathbb R^N\supset B^N(1).$

{\it Example.}   The composition of the toral Clifford embedding 
$\mathbb T^m\to B^{2m}(1)$ with the universal covering $\mathbb R^m\to \mathbb T^m$  followed   the Euclidean homothety $x\mapsto  (\sqrt n) x$ is an isometric immersion $\mathbb R^m\hookrightarrow B^{2m}(1)$ with curvature $\sqrt m$. Hence,
$$\mathcal D(m,2m)\leq\sqrt m$$ 

{\it Question.} Is  $\mathcal D(m,2m)$ equal to $\sqrt m?$ 
\vspace {1mm}

\textbf {1.1.A. Euclidean $\pmb{\mathcal D}\mathbf {(m,N)}$-Theorem.} 

$\bullet_{\geq 2m}$ If $N\geq 2m $, then
$${\mathcal D}(m,N)\leq\sqrt \frac {3m}{m+2}  +  C_o \frac{m} {\sqrt N},$$
where $C_o$ is a universal constant (see section 3). Moreover, if 
$N\geq 100m^2$, then 
$${\mathcal D}(m,N)=\sqrt \frac {3m}{m+2}. $$

$\bullet_{< 2m}$  If $m+1 \leq N<2m$, then 
$${\mathcal D}(m,N)\leq  6 \frac {m^\frac{3}{2}}{N-m} $$

{\it About the Proof.} The upper bound on $\mathcal D(m,N)$ is proven in section 3 for  $N\geq 2m $ and in section 4 for $ N\leq 2m.$

The lower bound on  $\mathcal D(m,N)$ and the issuing equality ${\mathcal D}(m,N)=\sqrt \frac {3m}{m+2} $ is proven in section 5 by reproducing  Petrunin's argument from [Pet2023]. 


\hspace {0mm} {\it Question.}  Is   ${\mathcal D}(m, m+1)$ bounded by $2m$?
 \vspace {1mm}
 
  \textbf { 1.1.B. $\mathbf\delta$-Approximation Corollary.} Let $X=X^m$ be a smooth manifold and $f: X\to \mathbb R^N$ a continuous map.
  
 $\bullet^\geq $ If $N\geq 2m-1$ then $f$ 
 can be $\delta$-approximated by smooth   immersions  
$$f^{}_\delta :X \hookrightarrow\mathbb R^N, \delta >0,$$   
with curvatures 
$$curv^\perp_{f^{}_\delta} (X) \leq \frac {1}{\delta} \left(\sqrt \frac {6m-2}{2m+1}  +  C_o \frac{m} {\sqrt N}\right)+o\left(\frac{1}{\delta}\right),\mbox { } \delta\to 0,$$ 
 where "$\delta$-approximated" means that 
 $$dist_{\mathbb R^N} (f^{}_\delta(x), f_0(x))\leq \delta,\mbox { }   x\in X.$$

  $\bullet^\leq $ If $X$ admits an immersion to $ \mathbb R^n$, $n< N$, and  $N\leq 2m$,
  then $f$ 
 can be $\delta$-approximated by smooth   immersions  
$$f^{}_\delta :X \hookrightarrow\mathbb R^N, \delta >0,$$   
with curvatures 
  $$curv^\perp_{f^{}_\delta} (X) \leq \frac {1}{\delta} \frac {6n^{\frac {3}{2}}}{N-n} +o\left(\frac{1}{\delta}\right).$$ 
  
 {\it Proof.}   Let $\phi: X=X^m\to \mathbb R^n$  be a  smooth  immersion \footnote{All  $X^m$ immerse to $\mathbb R^{2m-1}$,  if $m\geq 2$, 
  by the Whitney theorem.} and observe the following.
 
 \vspace {1mm}

\textbf {1.1.C. Stretching Lemma}. {\sf 
If $n\geq m+1$, then,  for all Riemannian metrics $g$ on $X$ and  all  positive  functions $\varepsilon (x)$, there exists an 
a $g$-expanding immersion $\psi: X\to \mathbb R^n$ {\it regularly homotopic} to $\phi$, i.e. it can be joined with $\phi$ by a $C^1$-continuous homotopy of  smooth immersion, and    such that 
$curv_\psi (X,x)\leq \varepsilon (x)$.}

{\it Proof.} If $X$ is compact,   scale   $\phi\to \psi=\lambda \phi $ and send $\lambda\to \infty$.

If $X$ is non-compact  and $n<m$ regularly homotop $\phi$ it to a {\it proper}  (infinity goes to infinity)  immersion with a use of Hirsch' immersion theorem    and  let  $\psi_\lambda:X\to \mathbb R^n$ be the composition of $\psi$ 
 with a    $\lambda(y)$-expanding map $:\mathbb  R^n\to \mathbb  R^n $, $y \in \mathbb  R^n$,   for a large and fast growing function $\lambda(y)$.\vspace {1mm}
 
Now, $\varepsilon$-approximate $f$ by a smooth map $f'_\varepsilon$ and add to it the  composed map of  $\delta^{-1} \psi_\lambda=\psi{\delta^{-1}\lambda}$ with an expanding map    $f_\odot:\mathbb R^n\to\mathbb R^N$ times $\delta$.
It is clear that if the function $\lambda(x)=\lambda_{f_\varepsilon}(x)$ is sufficiently large, depending  on  the norms of the fist and the second   differentials $||df'_\varepsilon(x)|| $ and   $||d^2f'(x)||$,  
then the curvature of 
this sum 
$$f_{\delta,\lambda} (x)=f'_\varepsilon(x) +\delta\cdot f^\odot\circ \psi_{\delta^{-1}\lambda }( \delta^{-1}x)$$ 
is bounded by 
$$\frac {curv (f^\odot)}{\delta}+
o\left(\frac{1}{\delta}\right)$$  
and the proof follows  with $\varepsilon \to 0$.

{\it \textbf { Remark}}  \textbf I. If $f=0$, and  $X$ immerses to $\mathbb R^n$, then the   above delivers an immersion $f_1$  of
$X$ to the unit ball $B^{n+1}=B^{n+1}(1)$ with a bound on the curvature of $f_1$ depending only on 
the dimension $m$ of $X$, e.g. 
$$min.curv^\perp(X,B^N)\leq \sqrt\frac {3(2m-1) }{2m+1}= \sqrt {3-\frac {6}{2m+1}   } \mbox  {  for } N\geq 100m^2. $$

Moreover,  we show in  section 3)the following.

\textbf {1.1.D.} {\sl As $N$ becomes very large {\it depending on the topology of $X$},
then  $min.curv^\perp(X,B^N)<\sqrt\frac {3(2m-1) }{2m+1}$.  In fact,
$$\lim_{N\to \infty} min.curv^\perp(X,B^N)\leq \sqrt\frac {3m}{m+2}=\sqrt {3-\frac {6}{m+1}   } .\leqno  { [\bf {N>>}]}$$}

\textbf {1.1.E.   Conjecture.}  {\sf If  $N\geq 100m^2$ then all $m$-manifolds $X$ admit immersions to the unit sphere 
 $S^N(1)$
   with curvatures $$curv^\perp(X\hookrightarrow S^N(1))\leq \sqrt {\frac {3m}{m+2}-1}=\sqrt\frac {2m-1}{m+2} .$$}

The  bound  $[\bf {N>>}]$, albeit unlikely,  may be optimal\footnote{Anton Petrunin [Pet2014]  proved it {\it is optimal}, see section 5.}   but our bounds on on $curv^\perp_{f_1}(X)$ for small $N$
are far from  optimal. For instance, 
  Clifford embeddings   of products of $l$ spheres to  the unit balls have   curvatures $l^\frac {1}{2}<<l^\frac {3}{2}$.

But the Clifford embeddings are not optimal either:  there   are  products of $l$ spheres, which admit codimension 1 (not $l$!) immersions with curvatures bounded by a universal constant, where the best available  -- we {\it don't know} if this is optimal -- such a constant is $1+2\sqrt \frac {3l-3}{l+1}$  according to the following.

\textbf {1.1.F. Codim 1 Theorem/Example}.(See section 4.2) Let 
$$X= S^{k} \times {\underset {l-1}{\underbrace {S^1\times...\times S^1}}}.$$
    If   $k \geq l^{l^4}$,\footnote{The hugeness of this  number is the product  of  my perfunctory  interpretation of  Hilbert's argument in [H1909]. }
then there exists an immersion 
$$F:X  \hookrightarrow B^{k+l}(1) $$ 
    with     
   $$curv_F^\perp(X)\leq 1+2\sqrt \frac {3l-3}{l+1}<4.5.$$

  \textbf {Remark II.} The   proof of the remark  \textbf I doesn't apply to immersions 
to  $\mathbb R^n$ without passing to  $\mathbb R^{n+1}$
 but this is taken care of by the following  
(see section 4.3).

\textbf {1.1.G.  Regular Homotopy/Approximation Theorem.} Let $f:X=X^m\to \mathbb R^n$ be an immersion. 
If $n>m$, then $f$ can be 
$\delta$-approximated by immersions $f_\delta: X\hookrightarrow\mathbb R^n$
which are regularly homotopic to $f$ and such that
 
 $$curv^\perp_{f^{}_\delta} (X) \leq \frac {500}{\delta} m^{\frac {3}{2}} +o\left(\frac{1}{\delta}\right).$$ 
 \vspace{2mm}

\textbf {1.1.H.  Remarks/Questions.}   We don't know how close  this inequality to the minimal values 
of the curvatures of codim1 immersions of products of spheres is. 

(a)  For instance  let $P^{l-1}$ be an $(l-1)$-dimensional manifold  diffeomorphic  to
a product of spheres where some of these have dimensions $\geq 2$. Then, if  $k>> l$, 
 there exist immersions 
$$F_\varepsilon :S^k\times P^{l -1} \hookrightarrow B^{k+l}(1) $$ 
    with     
   $$curv^\perp_{F_\varepsilon}(S^k\times P^{l-1})\leq 1+2\sqrt \frac {3l-3}{l+1}+\varepsilon $$
for all $\varepsilon>0$.

 But this  is {\it unclear} for $\varepsilon=0$, even  for
the product  $S^1\times S^{k}$, which   embeds to  the ball $B^{k+2}(1) $  with curvature 3   for  all $k$ and where   
  we {\it don't know}  if there are   immersions of $S^1\times S^{k+2}$ (or other closed  non-spherical manifolds of dimension $k+1$)  
 to  the unit ball $B^{k+2}(1) $   with curvatures $<3$.

(b) It is not impossible   according to what we know, that $m$-dimensional   products of spheres  of dimensions $\geq 2$ 
admit immersions to  $B^{m+1}(1)$ with curvature <100.

But the best we can do   (see section 4.1) are immersions with curvatures 
$ \lesssim m^\frac{4}{3}$.

\subsection {Equidimensional   Expanding Maps}

 \textbf {  Affine Expanding Maps.} The product of $r_i$-balls admits an {\it affine} equidimensional  
 expanding map  to the
$ R $-ball
$$f: \bigtimes_{i=1}^kB^{n_i}(r_i)\to B^N(R), \mbox { }N=\sum_in_i,$$
 {\it if and only if}  
 $$\sum_ir_i^2\leq  R^2,\leqno {\color {blue}[\sum r_i^2]}$$
 where  --  all this is, of course,  obvious -- in the case  of equality $\sum_ir_i^2=R^2$, such an $f$ is an {\it isometric embedding}.
 
 But -- this was pointed out to me by 
Roman Karasev-- it is  {\color{red!50!black}unlikely} that there is a simple criterion for the existence  of  
such  embeddings to cubes, not even for  
  rectangular solids, 
 $$\bigtimes_{i=1}^n B^1(r_i)= \bigtimes_{i=1}^n[-r_i,r_i]\to
[-\underline r ,\underline r]^n.$$\vspace {1mm}

\textbf{1.2.A.    Rolled  Band  Example.}  What is more interesting from our perspective is 
    {\it a $ (1-\varepsilon)$-expanding map, for a given $\varepsilon>0$,  from the infinite cylinder 
 $X= B^{n-1}(r)\times \mathbb R^1 $ to the ball $B^n(2r),$
 $$f_\varepsilon: B^{n-1}(r)\times \mathbb R^1  \to B^n(2r),$$}
where  this $f_\varepsilon$ comes as   the composition of  two maps.

  (1) The first map is the universal covering map from the  cylinder  $B^{n-1}(r-\varepsilon)\times \mathbb R^1 $
     to the {\it round  solid  torus} 
  embedded to the ball,
  $$f_1:B^{n-1}(r)\times \mathbb R^1 \to\mathbb T_{sld}(r, r-\varepsilon)\subset  B^n(2r),$$
  where this torus is equal to the   $(r-\varepsilon)$-{\it neighbourhood} of  a {\it planar  circle }
  $$ S^1(r)\subset B^n(2r)$$
  of {\it radius}  $r$, 
where  the center  of $S^1(r+\varepsilon)$ is positioned at   the center of the ball $ B^n(2r)$. 
  
  Observe that 
   the map  $f_1$ is {\it isometric} on the $(n-1)$-balls 
  $$B^{n-1}(r-\varepsilon)\times t \subset B^{n-1}(r-\varepsilon)\times \mathbb R^1,  \mbox {  $t\in \mathbb  R^1$}.$$
 
 (2) The second map  $f_2$ is the linear (scaling)  diffeomorphism 
$$f_2: B^{n-1}(r)\times \mathbb R^1\to B^{n-1}(r-\varepsilon)\times \mathbb R^1 \mbox { for }
 f_2:(s, t) \mapsto  \left (\frac{s}{1-\varepsilon}, \varepsilon ^{-1}t\right);$$
 where, clearly,   the composition
   $$B^{n-1}(r)\times \mathbb R^1\overset {f_2}\to B^{n-1}(r-\varepsilon)\times \mathbb R^1\overset {f_1}  \to \mathbb T_{sld}(r, r-\varepsilon)\subset  B^n(2r)$$
 is  the required  $ (1-\varepsilon)$-expanding map    $B^{n-1}(r)\times \mathbb R^1\overset {f_\varepsilon}\to   B^n(2r).$\vspace{1mm}

 \textbf{1.2.B  $[f\times f]$-Corollary.}  {\sf The Cartesian powers   of 
 $$f_\varepsilon: [-r,+r]\times \mathbb R^1  \to B^2(2r)\subset \mathbb R^2$$
deliver expanding maps
 $$ B^m(r)\times \mathbb R^m  \subset [-r,+r]^m\times \mathbb R^m \to B^{2m}\left(1+\frac {1}{\sqrt m}\right )$$
 for all $m=1,2,...$ and 
 $r< \frac {1}{\sqrt m}.$}\vspace {1mm}
 
  { \textbf {1.2.C.    $ \mathbf {\frac {1}{2}}$-Exercise}. } Show that if  $\underline r\leq 2r$, then the  
 cylinder  $B^{n-1}(r)\times\mathbb R^1$ admits {\it no expanding map} $f$  to the ball  
$B^n(\underline r)$. 

{\it Hint.}  (i){\it The axes} -- the central line $0\times \mathbb R^1$  of the cylinder  -- 
must go by $f$ to the {\it concentric ball} $B^n(\underline r-r)\subset B^n(\underline r)$.
 
(ii) The longest straight segment  with respect to the $f$-induced
flat metric   between pairs of points 
 on this  axes must have length $>2\underline r-r$.

The above   1.2.B is  generalized in section 4.2 as follows.

\vspace {1mm}

 \textbf {1.2.D. Rolled  Band into Ball Theorem.} {\it If $M\geq 100m^2$, 
 and $$r<\frac{\sqrt{m+2}}{\sqrt  {3m}+\sqrt{m+2}} \left(>\frac {1}{3}\right),$$
  then  the product   
 $ B^M(r)\times \mathbb R^m $ 
admits an  equidimensional   expanding map  to the unit  ball, 
  $$F_r:B^M(r)\times \mathbb R^m\to B^{m+M}(1).$$}

   {\it Remark/Question.} If $m=1$, then, by  the above ${\frac {1}{2}}$-exercise,   the bound  $r<1/2$ is optimal, but it is  not clear
   for $m=2$. 
   
   Here   the above inequality  for $m=2$, which allows expanding maps  from $B^4(r)\times \mathbb R^2$ to the unit ball 
   $B^{m+M}(1)$, where the supremum of the possible $ r$ is  
    $$\sup r=\frac {2}{\sqrt 6 +2}-\varepsilon (\approx 0.45),$$
    is implemented with $M=4$ by means of the  normal exponential map for the   2-{\it subtorus in   Clifford torus $ \mathbb T^3\subset B^6(1)$, which is  is normal to the principal diagonal in}  
  $\mathbb T^3.$

    Similarly  the  normal exponential map for the Clifford torus $ \mathbb T^2\subset B^4(1)$ leads to such maps 
   $B^2(r)\times  \mathbb R^2\subset B^4(1)$  with 
  $$\sup  = \frac {1}{1+\sqrt 2}\approx 0.41 <0.45,$$
   while the best $B^1(r)\times  \mathbb R^2\subset B^3$,
    where  
 $$\sup r=\frac {1}{3}<0.41,  $$  
   is obtained with the normal exponential map for the  standard  round torus in $\mathbb R^3.$
    
   And the only known upper bound on $r$  is  for $M=1$:
       $$r\leq \frac {\pi}{2\sqrt {\lambda_1(B^3(1))}}=\frac {\pi}{2j_{1/2}}=\frac{1} {2}>\frac{2}{\sqrt  {6}+2}\approx 0.45,$$ 
   where  this $\lambda_1$ is the first Dirichlet eigenvalue of the Laplacian in the unit 3-ball,
   and $j_{1/2}=\pi$ is  the first Bessel function zero(see section 5.1). \vspace{1mm}
   
{\it None of these four  inequalities is known to be} ({\it or not to be})  {\it optimal.}
\vspace {1mm}



\subsection { Remarks,  Acknowledgements and  the Plan of the Paper } 

The  lower bounds on curvatures
  of tori (see section 1.3) in concert with the  "natural symmetry"
 of Clifford's manifols
may lead one to believe   that such bounds  
persist in all codimensions. But when I mentioned   this to Fedia Bogomolov,
 "everything is possible in large dimensions" -- he responded.

Then my attempts to prove  lower bounds  on the curvatures of $m$-tori in $n$-dimensional balls for 
 $n\sim 2m$
 were arrested by what 
 Gilles Pisier explained  to me about  norms of generic  linear  families of 
selfadjoint operators.
  
Also Gilles pointed out to me on the criticality of dimensions $N\sim m^2$ (example 3.1 in [FLM1977]) 
  and the present state of art with Dvoretzky-Milman  inequalities  for the   
$l_p$-spaces was explained to me by  Grigoris Paouris  who also 
suggested to me  the relevance [K1995] for evaluation of the Kolmogorov diameter $D$ .

Then Bo'az Klartag and Noga Alon patiently  explained me  the  essential properties on the spherical designs and 
  construction of these  based on binary codes, allowing sharp bound on $D$ in moderately high dimensions. 
 We  present all this in section 2.

 In section 3, we show how bounds on the Kolmogorov $m$-diameter of the space $l^N_4$
translate to  corresponding inequalities for  curvatures $curv^\perp(X\hookrightarrow \mathbb R^{2N})$ for 
  submanifolds  $X$ in the Clifford tori 
$\mathbb T^N\subset \mathbb  R^{2N}$.

In section 4.1 we elaborate on the round torus construction from section  1  needed 
 for immersions below $4m-2$. 

In section 4.2. we exhibit   codim 1 immersions with small curvatures as boundaries of "tubular neighbourhoods" 
of immersion with  high codimension  constructed in the previous sections and   similarly construct   expanding 
maps in the cases indicated in section 1.2.

In section 4.3 we describe a  twisting procedure of immersed manifolds by regular homotopies  
with controlled curvature  and in section 4.4. we outline a similar procedure based on  
{\it Poenaru-Eliashberg's folding idea. }

\vspace{1mm}

 In section 5 we collect (mostly) known  bounds  on expansion and on the  curvature of immersions, including the recent  
 sharp  $\sqrt 3$-inequality by Petrunin.
 
 In section 6 we discuss  curvature problems similar to but  different from the 
 ones we address in the main body of the paper.


\section {Kolmogorov's $D=D(m,N,p)$, Hilbert's Theorem and Spherical Designs }


 

  \textbf {K-Diameter $ \mathbf {\sqrt [p]  {\mathbf {D(m,N,p)}}}$. } Let $||y||_{L_p}$, $y=(y_1,...,y_N)\in \mathbb R^N$  denote the  normalized 
norm   $l_p$,
  $$||y||_{L_p}=\left (\frac {1}{N}\sum_{i=1}^N|y_i|^p\right)^\frac {1}{p}$$

  Let $D(m,N,p)$ denotes the infimum of the numbers $D>0$ such that $\mathbb R^N$ contains an
   $m$-dimensional linear subspace $X$, such that
  $$   ||x||^p_{L_p}\leq D ||x||^ {p}_{L_2}, \mbox {  for all }  x\in X.$$
Observe that $D(1,N,p)=1$, $ D(m,m,p) = m^{\frac {p}{2}-1}$,  that  
 $D(m,N, p)$ is monotone increasing in $m$ and decreasing in $N$ and let 
   $$D(m,p)= D(m,\infty, p)=\lim_ {N\to \infty}D(m,N,p).$$  

  \textbf {2.1.A. Gamma Function Design Formula.} If  $p=4,6, 8...$, then  
 a simple $O(m)$-averaging argument,  shows that
$$D(m, p)=\frac {\int_{S^{m-1}} |l(s)|^pds}{\left( \int_{S^{m-1}} |l(s)|^2ds \right) ^\frac {p}{2}}
=\frac {m^{\frac {p}{2}-1} \cdot  3\cdot 5\cdot\cdot\cdot (p-1)}{(m+2)\cdot(m+4)\cdot\cdot\cdot(m+p-2)}, \leqno {[\Gamma/\Gamma]}$$ 
where $l(s)$ is a non-zero linear function on on the sphere. 

\textbf {2.1.B. Hilbert Connection}.  In his proof of the Waring problem,  Hilbert  shows the existence  of 
 {\it $M=\binom{m+p-1}{m-1}+1$ rational points $s_i\in S^{m-1}$} 
 and of  positive  {\it rational} weight $w_i>0$, $\sum_1^Mw_i=1$, such that  $\sum_iw_il^d(s_i)= \int_{S^{m-1}} l^d(s)d$
 {\it for all} linear functions on he sphere.
 
 This,  after partitioning  each  $s_i$  into $\Delta $ atoms for $\Delta $ being the  smallest common denominator $\cal N$
  of $ w_i$, becomes    what is no-a-days called  {\it spherical design}  of cardinality $N=\mathcal NM$  
  of $ w_i$,  which  yields (this is nearly obvious, see \textbf {2.1.C} below) the following.
        
       {\it  \textbf {$\mathbf {D(m,N)}$-Stabilization}}\textbf: 
      $D(m,N,p)=D(m,\infty,p)$ for all sufficiently large $N\geq =N_{Hilb} (m,p)(\leq\mathcal NM)$,  where -- to be safe let it be 
        rough--  $N_{Hilb}\leq m^{m^p}$.

{\it \textbf{Design Rationality}}\textbf:   If  $N\geq N_{Hilb}$ then   the space $l^N_p$ contains a {\it rational} linear subspace $X$ of dimension $m$, such that  
 $$||x||_{L_p}^p=D(m, p)||x||_{L_2}^{p} \mbox  {  for all $x\in X$.}$$

 \textbf {2.1.C. Spherical Designs and the Equality $D(m,N )=D (m,\infty)$ }
  
A  {\it design of even  degree $p=2,4,...$ and cardinality $N$} on the sphere $S^{m-1}$  is a map from a set  $\Sigma$  of cardinality  $N$ to the sphere, 
   written as $\sigma \mapsto s(\sigma)$, such that the linear functions $l(s)$ on the sphere  $S^{m-1}\subset  \mathbb R^m$
   satisfy
  $$\frac {1}{N}\sum_{\sigma\in \Sigma}  l^d(s(\sigma))=\int_{S^{m-1}} l^d(s)ds,  \mbox { } d=2,...,p,$$ 
  where $ds$ is the {\it $O(m)$ invariant probability} measure on the sphere.
  
Hence, the linear map from  the  space $\mathbb R^{m\perp}(=\mathbb R^m)$ of   linear 
functions  on the sphere  $S^{m-1}\subset  \mathbb R^m$  to
$\mathbb R^N=\mathbb R^\Sigma$ preserves both, the $L_2$ and the $L_p$-norms  and, by the above  
$ {[\Gamma/\Gamma]}$,

{\it the existence  a design of cardinality $N$ implies that 
 $D(m,N,p)=D(m,p)$}.\footnote {See [BB2009], [LW1993]  for more about it.}   \vspace {1mm}
  
 Non-rational designs, at least for $p=4$, are known to exit for $N<< N_{Hilb}.$ \vspace {1mm}
    
\textbf {2.1.D $2m^2$-Design Construction.} {\it If $p=4$, and if $m$ is a power of 2, then  there exists a  spherical designs of cardinality  $N=2m^2+4m$.}
 \footnote {This was stated and  proved in a  written  message by Bo'az Klartag to me.  Also, Bo'az pointed out to me that 
   the  Kerdock code used in   [K1995] yields  designs for $m=4^k$ and $N=\frac {m(m+2)}{2}$.}
  
  This, now for all $m$, shows that 
   $$ D(m,N,4)= \frac{3m}{m+2}
   \mbox  {  for $N\geq 8(m^2+m)$}.\leqno { \textbf {(i)}} $$
  
   [$\mathbb R^2$ in $l^3_4$]-{\it Example.} $D(2,N,4)=\frac {3}{2}$ 
 for $N\geq 3$, with  four (rational) planes $X\subset \mathbb R^3=l_4^3$, where $||x|||_{L_4}^4=\frac  {3}{2} ||x|||_{L_2}^4$:  these are 
  the normals to the vectors (1,1,1,), (1,1, -1), (1,-1,1), (1,-1,-1).
  \vspace {1mm}
 
  \textbf {2.1.E. $D(m,N)$-Inequalities.} If $N\lesssim m^2$, then upper bounds on $D^4(m,N,4) $ follow from 
   the corresponding 
   estimates in the 
  randomization proofs of the  Dvoretzky theorem for the $l_p$-spaces, where the   following inequality follow 
   from (the argument in) [PVZ2017]  as it was  spelled out in  details in a mesage  by Grigoris Paouris  to me.

 \textbf {(ii)} $D(m,N,4)  \leq 3 +const_{\mathbf {(ii)}}\frac {m^2}{N}$ for  $N\geq m^2$;\footnote {This follows from \textbf {(i)}
for   $N\geq 8(m^2+m)$  and, if $const_1$ is large,  also for  (some) $N\leq 8(m^2+m)$.  Besides,  the inequality 
 $D^4(m,m^2,4)\leq  const$ follows from (the proof of) example 3.1 in [FLM1977].}

\textbf {(iii) } $D(m,N,4)  \leq const_{\mathbf {(iii) }}\frac  {m^2 }{N}$  for 
 $2m\leq N\leq m^2$.\footnote { Since $D(m, N,4) \leq D(m, m,4) =m$ for all $m$ and $N$, the significance of this inequality for 
 $N\sim m$ depends on the value of $const_2$.} \vspace {1mm}

 \textbf {2.1.F.      $D(m,N)$ Concentration Property}.  The  existence of $m$-subspaces $X\in l^N_4$ in  
[FLM1977] and [PVZ2017], such that 
$$||x||^4_{L_4}\leq D||x||^4_{L_2}, \mbox { } x\in X, \leqno  [D]$$
is derived from a {\it  lower  bound the measure} of those $m$-subspaces $X\subset \mathbb R^N$, where 
this inequality fails for some $x\in X$. 

In particular,  the argument  used in  [FLM1977] implies that
the measure $\mu_D$ of those    $X\subset \mathbb R^N$ with respect to the 
$O(N)$-invariant  probability measure in the Grassmanian $Gr_m(\mathbb R^N$ where 
$||x||^4_{L_4}\geq D||x||^4_{L_2},$ for some  $x\in X$ satisfies:

{\it If ,  
$$D> \frac {3m}{m+2}$$  
then
$$\mu_D\to 0, \mbox { for $N\to \infty.$}$$}

{\it \textbf {Nash Connection}.} Besides applications to   lower bounds on  curvatures of   immersions (see next section), 
 Hilbert's argument,  combined with a Nash-like twist, leads to    $C^2$-smooth isometric  Riemannian  immersions with (large) 
  {\it prescribed} curvatures and also to a solution of the {\it differential geometric Warning problem}:
  
   construction of  {\it isometric $C^1$-immersions   of manifolds with  symmetric differential  forms of 
  degrees   $d>2$}, (see     2.4 (B)(4) on p. 205 in  [Gr1986] and [Gr2017]).

 
 \section{  Equivariant  Immersions  $\mathbb R^m \to S^{2N-1}$  and Euclidean $\mathcal D{(m,N)}$-Theorem
  for $N\geq 4n$ }


\textbf {3.A.  Curvatures of the Clifford Tori.}  Let $$\mathbb T^N\subset  S^{2N-1}\subset B^{2N}(1)\subset (B^2(1))^N\subset \mathbb R^{2N}$$
be the Clifford torus and observe that the second quadratic form of this 
torus  in the the ambient Euclidean space $\mathbb R^{2N}\supset S^{2N-1} \supset \mathbb T^N$,  
   regarded as a quadratic form
 with values in the normal bundle, 
is
$$\mbox {II=$\sqrt N\sum _{i=1}^N \nu_idt_i^2,$}$$ 
where $t_i$ are the  cyclic coordinates on the torus and $\{\nu_i\in T^\perp(\mathbb T^N\subset \mathbb  R^{2N})\}$  is the corresponding orthonormal frame of {\it normal} vectors to $\mathbb T^N$.
\vspace{1mm}

This, in terms of the orthonormal {\it tangent}  frame $\{e_i =\frac{\partial}{\partial t_i}\in T(\mathbb T^N)\}$, means that 
$$\mbox { II$: e_i\otimes e_i\mapsto \sqrt N \nu_i$  and   II$: e_i\otimes e_j\mapsto 0$  for $i\neq j$. }$$ 

Thus, the curvature of  $\mathbb T^N$ in $B^N$ along a unit  tangent vector  $\bar x\in T(\mathbb T^N)$,$$\mbox {$\bar x=\sum_i x_ie_i$, where $\sum_ix_i^2=1$,}$$
is 
$$\mbox { $curv^\perp(\mathbb T^N, \bar x)=||$II$(\bar x\otimes\bar x)|| =||$II$(\sum_i x_ie_i\otimes \sum_i x_ie_i)||=$}$$
$$\mbox {||II$(\sum_{ij} x_ix_j(e_i\otimes e_i)||= \sqrt N{ ||\sum_ix_i^2\nu_i||}=\sqrt N\sqrt {\sum_ix_i^4}=\sqrt N\frac{\sqrt {\sum_ix_i^4}}{||\bar x||^2}= $}$$
where  
$||\bar x||^2= |\bar x||^2_{l_2}   =\sum_{i=1}^Nx_i^2$. 

Hence,
$$curv^\perp(\mathbb T^N, \bar x)=\left (\sqrt[4] N \frac {||\bar x||_{l_4}}{||\bar x||_{l_2}}\right)^2= 
\left ( \frac {||\bar x||_{L_4}}{||\bar x||_{L_2}}\right)^2,\leqno \mbox {({\huge $\star$})}$$
where, recall,  the  $L_p$-norms refer to the  finite probability spaces with $N$ equal atoms, 
 $$||\bar x||_{L_p}=\frac { ||\bar x||_{l_p}}{\sqrt [p]N}.$$

\textbf {3.B. Proof of the Euclidean $\pmb{\mathcal D}\mathbf {(m,N)}$-Theorem 1.1.A  for $N\geq 2m$.} The above {({\huge $\star$}) implies the existence of an equivariant  isometric   immersion from the Euclidean $m$-space to the Clifford $N$-torus,
$$f^{\odot}:\mathbb R^m\to \mathbb T^N\subset S^{2N}\subset \mathbb R^{2N}$$ 
with the relative curvature  $curv^{\mathbf e}_{\mathbf E}(f^{\odot})$ 
(for the Euclidean metrics ${\mathbf e} $  in $\mathbb R^m$ and  ${\mathbf E} $  in $\mathbb R^{2N}$) 
equal to 
$\sqrt {D(m,N)}=\sqrt {D(m, N,4)}.$

Hence, 
$$\pmb{\mathcal D}\mathbf {(m,N)}\leq \sqrt {D(m,M)}$$
for all $m$ and $N\geq 2M$; thus
 the above $D(m,N)$-inequalities \textbf{(i),(ii),(iii)} yield the corresponding $\pmb{\mathcal D}\mathbf {(m,N)}$ inequalities in 1.1.

In addition to that,  if the $l^N_4$-space contains a  
rational $m$-subspace $X$ with $\frac {||x||^4_{L_4}}{ ||x||^4_{L_2}}=D,$
then  $\mathbb T^N$  contains  an $m$-subtorus with the ambient Euclidean curvature $\sqrt D$.

 \textbf {3.C.  Proof of $ [\bf {N>>}]$   from 1.1.D.} Embed   $X^m\to   \mathbb R^{2m}\subset \mathbb R^N$, $N>2m$, 
 and apply orthogonal transformations $o\in O(M)$ to $X.$ 
Since $X$ is compact  (non-compact manifolds are irrelevant here)  the   $D(m,N)$-concentration  2.1.F  implies that  
there exist an $o_\varepsilon\in O(n), $ 
 such that all tangent vectors $\tau\in o(T(X))\subset \mathbb R^N$ satisfy
$$||\tau||^4_{L_4}\leq \left (\frac {3m}{m+2} +\varepsilon_N\right)    ||\tau||^4_{L_2},\mbox
 { where $N_\varepsilon \to \infty $  for $N\to \infty$ }.$$

Thus, arguing as earlier,  the $\lambda$-scaled manifold $X$ imbeds to the Clifford torus 
$\mathbb T^N\subset S^{2N-1}\subset  \mathbb R^{2N}$ with 
 $$curv^\perp(X\hookrightarrow \mathbb R^{2N})\leq \left (\sqrt{\frac {3m}{m+2} 
  +\varepsilon_N}+ \epsilon_\lambda\right),\mbox  { where  $ \epsilon_\lambda \to 0$  for $\lambda\to\infty$ }$$ 
 and the proof follows.

 \textbf {3.D. $\delta$-Approximation in  Non-Euclidean Riemannian Manifolds.}   The derivation  of the 
 $\delta$-approximation  from expanding Euclidean maps in section 1.1 easily generalizes, albeit with limitations, to Riemannian manifolds  as follows.

  \textbf {Theorem.} {\sf Let $Y$ be a complete Riemannian manifold\footnote{One may allow a boundary,  but this is a minor problem.} with the sectional curvature $|sect,curv^\perp(Y)|\leq \kappa^2$
 and let $f:X\to Y$ be a continuous  map. 
 
 If the induced bundle  
 $f^\ast(T(Y))\to X$ contains a  subbundle isomorphic to $X\times \mathbb R^N$, (i.e. a trivial one)  and if $X$ admits an immersion to $\mathbb R^N$, e.g.  $2m-1\leq N\leq dimY-dim (Y)-1$, then, 
 for all positive $\delta\leq \frac {1}{2\kappa}$,    the map $f$ can be $\delta$-approximated by   immersions $f_\delta:X\to Y$, such that 
  $$curv _{f_\delta}(X) \leq \frac {1+2\kappa}{\delta} \sqrt {\underline{\mathcal  D}_{}(m,N)},$$
  where
  $$\underline{\mathcal  D}(m,N)\leq \frac {3m}{m+2} +
   const \frac {m}{\sqrt N}\mbox {   for } N\geq 2m$$ 
  
 and $$\underline{\mathcal  D}(m,N)\leq \frac {6m^\frac {3}{2}}{N-m} \mbox {   for } N\leq 2m.$$}
 
 {\it Proof. } Proceed as in the proof of 1.1.B, where instead of adding $ \delta\cdot f^\odot\circ \psi_{\delta^{-1}\lambda }$ to  $f'_{\varepsilon}$ we  the compose    exponential map  
 with a (fiberwise injective)   bundle homomorphism  from the trivial bundle $X\times  \mathbb R^N$ to $X$ over the smooth map 
 $f'_{\varepsilon}$,  (this map $\varepsilon$-approximates $f$).

  \subsection  {Subtori in Non-Equilateral Clifford Tori} 
  All  invariant $N$-tori in the sphere $S^{2N-1}\subset\mathbb R^{2n}$   are  (equal, up to  isometries of $S^{2N-1}$, to)  the  orbits of the
   product action of  $N$-copies of the standard action of  $\mathbb T^1$ in the plane. where these orbits  
     are equal to the  non-equilateral Clifford tori 
   $$\mathbb T^N(\bar r) \bigtimes _{i=1}^NS^1(r_i), \mbox { for  $\bar r=(r_1,...,r_N)$, where $||\underline r||^2= 
    \sum_i r_i^2=1$}$$
     
  Then, similarly to the above {({\huge $\star$})}, the values of the curvature operator  of this torus at the unit tangent vectors $\bar x =(x_1,....x_ N) \in T(\mathbb T^N(\bar r))$ are 
$$curv^\perp(\mathbb T^N_{\bar r}, \bar x)= \left |\left |\sum_i \frac { x_i^2}{r_i} \nu_i\right |\right|  =\sqrt {\sum_i\frac {x_i^4}{r_i^2}} \leqno \mbox {({\huge $\star$}\hspace{-0.6mm}$\bar {\mathcal R})$}$$
where, 
if  all $r_i=\frac{1}{\sqrt N}$,  this reduces to  {({\huge $\star$})} for
$$\sqrt {\sum_i\frac {x_i^4}{r_i^2}}=  \sqrt \frac{\sum_i|x_i|^4}{N}$$
and where we denote
 $$||x||_{L_4(\bar r)}=\sqrt[4] {\sum_{i=1}^N\frac {x_i^4}{r_i^2}}$$

  \textbf {3.1.A.} {\it \textbf {Conclusion}.}
   {\sf There is a one-to-one correspondence between
 {\it $$\mbox {equivariant   $\mathbb R^m\subset S^{2N-1}$ with
{\color {black}$curv^\perp(\mathbb R^m)< \alpha$}}$$ and 
 pairs $( \bar r, X)$,   
   where
   $\bar r=(r_1,...,r_N)$ is a unit vector with positive entries, 
$$\sum_{i=1}^N r_i^2=1,\mbox { }r_i>0, $$
 and subspaces $X\subset Y=\mathbb R^N=l^N_2$ is a    

such that all  $x\in X$ satisfy
$$||x||_{L_4(\bar r)}< \sqrt \alpha\cdot||y||_{L_2}, $$}}
where, recall,  the $L_2$-norm of $y\in Y$, including $y\in X\subset Y$, 
is $$||y||_{L_2}=\sqrt\frac{ \sum_{i=1}^N y_i^2}{N}=\frac {||y||}{\sqrt N}. $$

 Conceivably,   $m$-torical orbits  not contained in  $\mathbb T^N_{Cl}$,
    e.g. those {\it maximizing the $m$-volumes}  of the respective $m$-tori actions,
  may have slightly smaller curvatures  than Kolmogorov's $D(m,N)$, that  is, as 
  we know, is equal to the infimum   
  of the curvatures of  $m$-subtori in $\mathbb T^M_{Cl}$.

  This can be stated with the  $\bar r$-counterpart of Kolmogorov's   $D(m,N)$,
   denoted $\largelozenge(m,N)(\leq D(m,N))$ that is  the infimum of the suprema of 
  the  ratios of the 
  two norms:
 $$\largelozenge(m,N)=\inf_{Y,\bar r}\sup_{0\neq y\in Y}  \frac {||y||_{L_4(\bar r)}}{||y||_{L_2}},$$
  where the infimum is taken over all $m$-dimensional  linear  subspaces   $Y\subset \mathbb R^N$ and all positive unit vectors $\bar r$.\footnote{Grigoris Paouris  has sent to me a message  with an evaluation  of $\largelozenge_{\bar r}=\inf_{Y}\sup_{0\neq y\in Y}  \frac {||y||_{L_4(\bar r)}}{||y||_{L_2}},$ 
  for several classes of  $\bar r$.}

 {\it   Question.} Is, ever,
  $\largelozenge(m,N)<D(m,N)$?\footnote{This can' t happen 
  for large $N>>m^2$ by Petrunin's inequality.}
  
  \vspace {1mm}

\textbf { The space $\mathbf {\mathcal I_\alpha=\mathcal I(m,N,\alpha)} $} {\it of isometric equivariant immersions 
 $\mathbb R^m\hookrightarrow S^{2N-1}$  with curvatures $\leq \alpha$}
 is a semi algebraic subset in the (Euclidean)  space
 $J_N(m,N)$  of $N$-jets at 
 $0\in\mathbb R^m$ of smooth maps $ \mathbb R^m\to \mathbb R^N$\footnote{The space $J_k(m,N)$ is isomorphic to the space of 
 polynomial maps $\mathbb R^m\to \mathbb R^N $ of degrees $\leq k$}, which
  is invariant under the action of the 
 orthogonal group $O(2N)$, and where
 the $O(2N)$-orbit of an $I\in \mathcal I$  in $S^{2N-1}$ is equal to

 {\sf $W_I\backslash O(2N)\slash \mathbb T^N$,
where $W_I$ is the subgroup of the Weyl group of $O(2N)$,  which preserves  $I$,  (this is empty for generic $I$).\footnote{ 
  The corresponding space $\mathcal X(m,N, \sqrt \alpha)$ of $m$-subspaces $X$ in $L^N_4$  with 
 $\frac {||x||^4_{L_4}}{ ||x||^4_{L_2}}=\sqrt \alpha,$ which, albeit being also  semi algebraic,    has more combinatorial flavour than $\mathcal I$.}}

There can be something geometrically interesting in the     $O(N)$-topology of  $\mathcal I_\alpha$ depending on  $\alpha$, but all one can say off hand  is the 
{\it Petrovsky}-(Thom-Milnor)  {\it bound on    the homology of }$\mathcal I_\alpha$   by the algebraic degree of this set.


\section  {Normal Immersions in Small Codimensions}

\subsection  {Proof of Euclidean $\mathcal D(m,N)$-Theorem for  $N\leq 2m$}



\textbf {{\Large$\rtimes $}-Construction.}  Let      $\phi_1:X_1=X_1^{ m_1}\hookrightarrow  \mathbb R^{m_1+n_1}$, be an immersion with a {\it trivial} normal  normal bundle,
where this "triviality" is implemented by a smooth map 
$$\Phi_1: X^1\times \mathbb R ^{n_1}\to \mathbb R^{m_1+n_1}$$ and let 
$\phi_2: X_2=X^{m_2}\to \mathbb R ^{n_1}$ be another immersion.
If $\phi_2$ lands in the $r$-ball in $  \mathbb R ^{n_1}$  for some  $r_1>0$,
$$\phi_2(X_2)\subset B^{n_1}_0(r)\subset \mathbb R ^{n_1}$$ and 
$$curv^\perp_{\phi_1}(X_1)\leq \alpha_1<\frac {1}{r},$$
then  the composed map  $(x_1,x_2)\mapsto \Phi_1(x_1, \phi_2(x_2)$ is an  {\it immersion,} 
say 
$$\phi_1\rtimes \phi_2: X_1\times X_2\to    \mathbb R^{m_1+n_1}.$$
 
Recall that the {\it normal connection $\nabla^\perp$ } in the (trivial) normal bundle 
$$X_1\times  \mathbb R^{n_1}=T^\perp(X_1)=T( \mathbb 
R^{m_1+n_1})\ominus T(X^1)\to X_1$$
is defined by the field  $\tau^\perp$of tangent $m_1$-planes in $X_1\times \mathbb  R^{n_1}$,
 which are normal  to the 
Euclidean fibers $x_1\times \mathbb R^{n_1}$ with respect to the
 (flat) Riemannian metric induced by the map
$\Phi_1: X^1\times \mathbb R ^{n_1}\to \mathbb R^{m_1+n_1}$.
 
{\it \textbf  {Flat Split Bundles and {$\nabla^\perp$-Trivial Immersions}} } The connection  $\nabla^\perp$ is called {\it flat split} if the map $\Phi_1$ is  {\it $\nabla^\perp$-parallel}
that is the field  $\nabla^\perp$ is normal to the fibers $x_1\times \mathbb R^{n_1}$  with 
respect the product metric in  $X^1\times \mathbb R ^{n_1}$ and  the immersion $ \phi_1$ 
is called  {\it $\nabla^\perp$-trivial} in this  case.

 \textbf {4.1.A.  List of $\nabla^\perp$-Trivial  Examples.} (a) Immersions $\mathbb R^1\to \mathbb R^n$ are $\nabla^\perp$-trivial.

 (b) Codimension 1 immersion of orientable manifolds, $X^m\to \mathbb R^{m+1}$,
  are   $\nabla^\perp$-trivial.
 
 (c) Equivariant immersions of tori, $\mathbb T^m\to \mathbb R^n$, are  $\nabla^\perp$-trivial.

(d) Direct products of $\nabla^\perp$-trivial-immersions $\phi_i :X_i\to \mathbb R^{n_i}$ 
$$\bigtimes_i \phi_i :\bigtimes_i X_i\to \mathbb R^{\sum_in_i}$$
are $\nabla^\perp$-trivial.

(e) The above "semidirect products" 
$\phi_1\rtimes \phi_2: X_1\times X_2\to    \mathbb R^{m_1+n_1}$
of  $\nabla^\perp$-trivial  $\phi_1:X_1\to \mathbb R^{m_1+n_1}$ and
 $\phi_2:X_2\to \mathbb R^{n_1}$ are  $\nabla^\perp$-trivial.

 \textbf {4.1.B.} (Obvious) {\Large $\rtimes$}-\textbf{Normality  Lemma.} Let   $\phi_1:X_1\to \mathbb R^{n_1}$
  and $\phi_2:X_2\to \mathbb R^{n_1}$ be  $\nabla^\perp$-trivial  immersions. Then:
 
 $\bullet_{norm}$ If $\phi_1:X_1$
  and $\phi_2$  are normal (see 1.1) then $\phi_1\rtimes \phi_2$ is also normal.

$\bullet_{curv}$  If $\phi_2(X_2)\subset B^{n_1}(r)\subset  \mathbb R^{n_1}$,
 then
$$foc.rad_{\phi_1}\rtimes \phi_2(X_1\times X_2)\geq \min ( foc.rad_{\phi_2}(X  _2), foc.rad_{\phi_1}(X_1)-r)$$
 and  in the normal case  the  relative curvature  of  $\phi_1\rtimes \phi_2$    (as well as the curvature $curv (X)=foc.rad(X)^{-1}$ itself), satisfies 
 the corresponding inequality.
$$curv^\perp(\phi_1\rtimes \phi_2)\leq  (\min ( curv^\perp(\phi_2)^{-1}, curv^\perp(\phi_1)^{-1}-r))^{-1}. $$

\vspace {1mm}

\textbf {4.1.C. Torus-by-Torus Construction.}  Let
$$[-1,1]\times \mathbb T^1\to [-2,2]^2\supset B^2(2)$$
be the map  obtained by  rotation of the segment $[0,2$ around the origin in the plane (which is an immersion away 
from the "interior" boundary circle)   and let
$$f_1=f_0^{\times k}: [-1,1]^k\times \mathbb T^k=([-1,1]^k\times \mathbb T^1)^k= \to   ([-2,2]^{2})^k=  [-2,2]^{2k},$$
$$f_2:[-1,1]^k\times \mathbb T^{3k}=[-1,1]^k\times \mathbb T^{k}\times \mathbb T^{2k}\to    [-2,2]^{2k}\times \mathbb T^{2k}
 = ([-2,2]^k\times \mathbb T^k)^2\to  [-4,4]^{4k}$$
  
  \hspace{10mm}.............................................................................................................

$$f_i:[-1, 1]^k\times \mathbb T^{k2^i-k}\to [-2^i,2^i]^{k2^i}.$$

 It follows by the construction, that this map is  normal  and that the normal  exponential map of the central torus 
 $$ \mathbb T^{2^i-1}=0\times \mathbb T^{2^i-1}$$ (immersed actually embedded) to
 the cube   $  [-2^i,2^i]^{k2^i}$
 is injective in the interior of  $[-1, 1]^k\times  \mathbb T^{k2^i-k}.$
 Hence, the curvature  of this torus and the (relative) curvature of the immersion $f_i$ are bounded by $1$ and the 
 corresponding scaled 
 map 
 $f: \mathbb T^{k2^i-k} \to B^{k2^k}$  satisfies
 $$curv_F^\perp( \mathbb T^{k2^i-k}) =curv^{\mathbb T^{k2^i-k}} (f)   \leq  2^i\cdot  \sqrt {k2^i},$$
 or, in terms of $m={k2^i-k}$,  
 $$curv_F^\perp( \mathbb T^m)\leq \left(\frac {m}{k}+1\right) \sqrt{m+k}, $$
 which implies for all $m$  and $k\leq m$:
 $$curv_F^\perp( \mathbb T^m)=curv^{\mathbb T^m}(f)<6\frac {m^\frac {3}{2}}{k}.$$
  The proof of theorem 1.1.B is concluded.
 \subsection {Proofs of the Codim 1  and the Rolled Band 
 Theorems}.
 
 
  Let  $f:X^m\to Y$ be an immersion with  $foc.rad_f(X) = R$ and  $S^\perp(r)(X)\to X$ be the bundle of normal  $r$-spheres   $ S_x^{N-m-1}(r)\subset T_x^\perp(X)= T_{f(x)}(Y)\ominus T_x(X) =\mathbb R^{N-m}.$

If $r<R$  then the normal exponential map 
 $ E:S^\perp(r)(X)\to Y$ is an immersion, where
 $foc.rad_E(S^\perp(r)(X))=\min (r,R-r).$

For instance, if $X\to B^N(1)$ is an immersion with trivial normal bundle and $curv_F^\perp(X)\leq $, then
the immersion
 $$E_f: \left(1+\frac {1}{2c}\right)^{-1} E:X\times S^{N-m-1}=S^\perp\left(\frac {1}{2c}\right)(X)\to B^N(1)$$
 has 
 $$curv^\perp_{E_f}\left(X\times S^{N-m-1}\hookrightarrow B^{N-m-1}\right)\leq
 {2c}\left(1+\frac {1}{2c}\right)=(2c+1).$$
 
 \textbf {4.2.A. Codim1 Conclusion.} This, applied to immersions of tori $\mathbb T^{l-1}\to B^N(1)$
 with large $N$ curvature  
 $\mathbb T^{l-1}=\sqrt \frac {3 (l-1)}{l+1}$, yields  codimension codimension one immersions with small curvature as stated in 1.1.G.
 
 \textbf {4.2.B.  Generalization from $l$-Tori to $l$-Polyhedra. } Given a compact polyhedral (or cellular) space $P$ of dimension $l$, there exists a compact $N$-manifold $X$, for all $N\geq 2l-1$, such that:  
  
 $\bullet_{P}$ there is a continuous map $K\to X$, which is a homotopy equivalence in dimensions $<N/2$, i.e.
  this map induces isomorphisms of the  homotopy groups,
  $\pi_i(P)\to \pi(X)$ for $i<N/2$;

 $\bullet_{200}$ if $N\geq 200 l^2$ then, for all $ \varepsilon>0$, 
 $X$ admits an immersion to  $ B^{N+1}(1)$ with
$$curv^\perp(X\hookrightarrow B^{N+1}(1))\leq 1+\sqrt \frac {3l}{l+2} +\varepsilon.$$
In fact, the boundary  of the regular neighbourhood of $P$ embedded to  $\mathbb R^{N+1}$ can be taken for $X$.
 
{\it  Embedding Remark. }  This, $X$, by its very  construction, {\it embeds} to 
 $\mathbb R^{N+1}$, but  one can show (section 5.3)
  that there is {\it no universal bound  on the curvature} of  embeddings of $X$ to the unit ball in $\mathbb R^{N+1}$.
  
  For instance if $P$ is a 
 connected sum of different lens  spaces, e.g. 
$$P_k=\mbox{ {\large\#}}_{i=1}^kS^3/\mathbb Z_{p_i},$$
where $p_1< ...<p_i<...<p_k$ are prime numbers,
 then the curvatures of all smooth embeddings $F:X\to B^{N+1}(1)$ satisfy:
 $$curv_F^\perp(X)\geq \log\log(k)/N^N.$$
 
 {\it Question.}  What, roughly, is the minimum of the  curvatures 
 of embeddings $\mathbb T^l\times S^{N}\to B^{N+l+1}(1)$?
 (See section 6.3 for more about it.)
 \vspace{2mm}

\textbf {4.2.C. The proof of the "rolled band  theorem} proceeds similarly to the above.

 Let  $f:\mathbb  R^m\to B^{m+M}(1) $ be an  immersion   with curvature bounded 
by $\mathcal D=\mathcal D(m,m+M)$ as in 1.1. 
let 
$$e=e_f:\mathbb R^m\times B^M(r) \to \mathbb
  R^m\to B^{m+M}(1+r), \mbox {  } r<\frac {1}{\mathcal D}, $$
    be the normal exponential map for 
  $\mathbb R^m$ immersed to $\mathbb R^{m+M}\supset B^{m+M}$ and let 
  $$E_\lambda : \mathbb R^m\times B^M(r) \to \mathbb
  R^m\to B^{m+M}(1)\mbox { for } (x,b)\mapsto  (1+r)^{-1} e(\lambda x,b).$$

If $\lambda$ is sufficiently large, then the map $E_\lambda$ is expanding in the $\mathbb 
R^m$ directions, i.e. it  expands $\mathbb R^m\times b$ for all $b\in B$   and since it is isometric in the $B^M$-directions it is expanding on $\mathbb R^m\times B^M(r)$... except for one problem: 

{\sf the normal $M$-ball  bundle $B^\perp(r)\to \mathbb R^m$ 
 of the  immersed $\mathbb R^m\hookrightarrow \mathbb R^{m+M}$ is trivial, it is  indeed, isomorphic to the product  $\mathbb R^m\times B^M(r)$ but the map $(x,b)\mapsto 
 \lambda (x,b)$ is not necessarily expanding  with respect to the 
(Euclidean) metric induced by the exponential map. (Look at the planar  map
  $(x,y)\mapsto(0, 10x+y)$  }

Fortunately, the normal bundles of our immersions constructed in sections  and 3. are {\it flat split}, (see 4.1)  the map $E_\lambda$ is expanding and it can be taken for the required $F_r$ in.

 \textbf {4.2.D.   Expanding Maps $F_r$  for all $m$ and $M$}. The above argument delivers expanding maps
 $F_r:\mathbb R^m\times B^M(r)\to B^{M+m}(1)$ provided 
   $r\leq \left ( 1+{\Delta}\right)^{-1}$, 
   where $\Delta $  is taken according to the $\mathcal D(m,N)$ inequalities (see  section 1.1 and  3).

$$\Delta =\sqrt \frac {3m}{m+2}  +  C_o \frac{m} {\sqrt M}, \mbox { for $M\geq m $},$$
and

$$\Delta=  6 \frac {m^\frac{3}{2}}{M} \mbox { for $M<m$}.$$
.

 \subsection {Proof of the  Regular Homotopy/Approximation Theorem.}  
 \textbf {Step 1. {\it Slicing.}} Given an immersed manifold  
 $$X=X^m \overset {\phi}\hookrightarrow \mathbb R^n,\mbox { } n>m,$$, and (small) positive numbers $\varepsilon,\delta>0$ there exists an 
  immersion   $$X \overset {\varphi}\hookrightarrow \mathbb R^n$$ 
  regularly homotopic tp $\phi$, such that 
 
 $\bullet_{curv_\varphi}$   $curv_\varphi(X)\leq \varepsilon$,

  $\bullet_\delta$ the first coordinate function $y_1(x)=y_1(\varphi(x))$  of
 $y=\phi(x)\in \mathbb R^n=\{y_1,....y_n\}$ 
is  proper Morse,   where 
   there are   no critical points of $y_1$  on the $\delta i$  levels  of $y_1$  for integer $i =...-2, -1,0,1,2...$,  i.e.
the hyperplanes  where $y_1=\delta i$ in $\mathbb R^n$  are transversal to   $ \varpi(X)\subset \mathbb R^n$ and

 $\bullet_\varepsilon$ the curvatures of these $ \delta i$ levels are bounded by $\varepsilon$.
 \vspace {1mm}
  
 {\it Proof.} If $X$ is compact, then  $\bullet_{curv_\varphi}$  achieved  achieved by scaling: 
 $x\mapsto \lambda \phi(x)$ for a large $\lambda$ and then one 
 gets $\bullet_\delta$  by a preliminary  generic rotation of 
$\phi(X) in \mathbb R^n$,  where
   then the critical values of $y_1(x)$ moved to 
the  centers  of the  segments  $[\delta i,\delta (i+1)]$,  
 let  
 $\frac {1}{\delta}=o(\lambda)$and   conclude the 
 proof  with the following obvious (but essential)

 \textbf { 4.3.A.  Levels Curvature Sublemma}. Let $y(x)$ be a Morse function on a compact Riemannian manifold $X$ and $x_0$ be a critical  point, where $y(x_0)=0$. Then the 
curvatures of the $\delta$-levels $f^{-1}(\delta)\subset X$ satisfy
$$curv^\perp(f^{-1}(\delta)=o\left (\frac {1}{\delta}\right).$$

\textbf {Step 2.  {\it Zigzag Folding and Compression.}}  Reflect  the $X$-bands $y_1^{-1}[\delta i,  \delta( i+1)\subset X$    in the hyperplanes 
$y_1=\delta i$, $i\in\mathbb Z$,  and thus "compress" $\varphi(X)$ to a zigzag map $\zeta$ from $X$ to the Euclidean $\delta$-band between a pair of such hyperplane, say between  $y_1=0$ and $y_1=\delta$. 
 
 \textbf {Step 3. {\it Twisted  Regularization  with Controlled Curvature.} }
 There  exists  a smooth $10\delta$-approximation of $\zeta$ by a smooth immersion $\zeta_\circ: X\to \mathbb R^n$, such that  
 
 $\bullet_\epsilon$  the immersion $\zeta_\circ$ is equal to $\zeta$ outside the $\epsilon$-neighbourhood of the{\it corners of 
 $\zeta$}, that is the subset $y^{-1}_1(\delta \mathbb Z)\subset X$,
 where $\epsilon>0$ en is a given number which may be taken much smaller than $\delta$;
 
 $\bullet _{reg}$ the immersion $\zeta_\circ$ is regularly homotopic to $\varphi$,
 
 $\bullet _{curv/ \delta}$ the curvature $\zeta_\circ$is bounded by 
 $\frac {1}{\delta}$
 
 {\it Proof}. To see how it works, let $\theta_\circ$ and $\theta_\neswbipropto$  be two immersions of  the circle to the plane, 
 each having a single corner point, both with the same corner  angle.  If we align these corners  properly and attach the immersions one to another at the corner  points, we  obtain a composed smooth  immersion $\theta_\ast$  where, if  $\theta_\neswbipropto$ is $\neswbipropto$-shaped, this
  $f_\ast$ is  regularly homotopic to $f_\circ$.
 
Now, in he case of a corner along a hypersurface $X_i =\varphi^{-1}\delta)$ attach the product $X_i\times  \neswbipropto$ to $\zeta (X)$ along this corner  and  by doing it to all $X_i$ we obtain a smooth immersion  regularly homotopic to$\varphi$ where the conditions $\bullet_\epsilon$ and  $\bullet _{curv/ \delta}$  are easily achievable $10\delta$ close 
to $\zeta$. Details are  left to the reader.

\textbf {Step 4. {\it Rolling  Bands into Balls.}}  The  band $\mathbb R^{n-1}\times [-10\delta,11\delta]\supset \zeta_\circ(X) $ is mapped to $B^n(1)$ by "rolled band" immersion $F_r:\mathbb R^n\times [-r,r] \to B^n(r)$ for $r$  from 1.2.D, where $F_r$ is restricted to the sub-band $\mathbb R^n\times [-r/2,r/2]\mathbb R^n\times [-r,r]$ and where we then  let
$\delta=\frac {r}{42}.$

In order estimate the curvature of the composed map 
$\Phi=F_r\circ\zeta_\circ$, 
$$X\overset {\zeta_\circ}\to \mathbb R^n\times [-r/2,r/2]\overset {F_r}\to B^n(1),$$  
by  $curv^\perp_{\zeta_\circ}(X)\leq c= \frac {1}{\delta}$
we recall the construction of the underlying normal immersion
$$f=F_r|_{\mathbb R^{n-1}\times \{0\}} : \mathbb R^{n-1}\to B^n(1-r),$$where 
$curv (f) \leq 6(1-r)^{-1} n^\frac {3}{2}$ and where also (the differential of) this map
has {\it controllably  bounded anisotropy},
$$\frac {||d\tau_1||}{||d \tau_2||}\leq 2n$$
for all unit tangent vectors $\tau_1,\tau_2\in T(X)$.
 It follows that the curvature $curv_\phi^\perp(X)$ 
is bounded essentially in the  same way as that of $F_r$, 
$$curv^\perp_{\Phi}(X) \leq420 n^\frac{3}{2},$$
and the corresponding approximation inequality  follows as  in the  proof in the genera case of the  $\delta$-approximation theorem. (This $\delta$ and that in $[-10\delta,11\delta]$,
 albeit similar, are not the same.)

  \textbf { 4.3.B.  Immersions to non-Euclidean $Y$}. The above  argument, unlike the proof of the 
 the $\delta$-approximation theorem as explained  in 3.D doesn't  generalize to immersions from $X$  to  general Riemannian manifolds $Y$.
 
  Yet,  a combination of the above "twisted  regularization"  on the top of  a routine  induction by  skeleta   delivers the following.

  \textbf {4.3.C. Rough Exponential Bound on Curvature}. {\sf Let $Y$ be a complete Riemannian manifold with 
$|sect.curv^\perp|\leq \kappa^2$ and let
$f:X=X^m \hookrightarrow Y$ be a smooth immersion.}

{\it  If $dim(Y)>m$  then, for all positive $\delta\leq \frac{1}{\kappa}$, the map $f$ can be $\delta$-approximated by immersions $f_\delta : X\to Y$,
 which are regularly homotopic to $f$ and such that
 $$curv^\perp_{f_\delta}(X)\leq\frac {(1+\kappa) 100^m}{\delta}.$$}

   \subsection { Unfolding Folds  and other Singularities.} 
 
 Below is another proof of  the regular homotopy/approximation theorem for {\it orientable hypersurfaces}, which leads to a better,  possibly sharp in some cases, bounds on the curvature.

\textbf {Unfolding Lemma.}
{\sf Let $X=X^m$ be an orientable manifold and
 $f:X\to \mathbb R^{m+1} $ be an immersion. Then, for all $\varepsilon>0$, there is an immersion,  
 $$\zeta_\circ:X\to \mathbb R^m\times [-1, 1],$$  
 which is regularly homotopic to $f$ and such that
  $$curv^\perp_{\zeta_\circ}(X)\leq 1+\varepsilon.$$}

 {\it Proof.} 
Apply Poenaru's $h$-principle for pleated maps (see (C) on p.56 in  [Gr1986]),  and obtain  a smooth map  $f_1:X\to \mathbb  R^{m+1}$ regularly homotopic to $f$, such that the only singularity of the normal projection $\zeta :X \to \mathbb R^m\subset \mathbb R^{m+1}$
 is a folding along  a smooth hypersurface $\Sigma=\Sigma^{m-1}\subset X$.

Make the curvature of the immersion $\zeta: \Sigma\hookrightarrow \mathbb R^m $ as  small as you wish by $\lambda$-scaling  as we did earlier  and thus also separate  different part of $\Sigma$ far one from another, such that,  on the balls of large radii $R\sim\lambda$ in $X$, the scaled  map is $\varepsilon$-close to the  standard fold 
$(x_1,...,x_m)\mapsto (x_1,...,x^2_m).$

 "Unfold" $\zeta\leadsto \zeta_\circ =(\lambda \zeta, y)\in \mathbb R^{m+1}$,
where $y:X\to \mathbb R$ is a smooth function on $X$, which, in the obvious normal coordinates, depends only on the last coordinate  $x=x_m$, where it is $\varepsilon$-close to a lift $ \eta_\circ:\mathbb R \to \mathbb R_+\times [-1,1]$ of the standard fold 
 $\mathbb R\to\mathbb R_+$,  $x\mapsto  y=x^2$, 
 where    $\eta_\circ(x)=(x,y(x))$ and where

 {\sl  the $x$-segment  
 $[-1,1]$ is sent by $\eta_\circ$ to the semicircle in the 
 half plane  $\{x,y\}_{y\geq 0}$ and 
  $\eta_\circ(x)=-1$  for $x< -1$ and   $\eta_\circ(x)=1$ for $x >1$.}

 Conclude the proof by rolling the band $\mathbb R^m\times [-1,1]$ into the ball as in the above   step 4.

 {\it Remarks.} (a) Our unfolding with controlled curvature quantifies 
  a single step in {\it removal of the singularities}  argument
   (see [GE1971]  and  section 2.1 in [Gr1986].)

To do the same for all step  and thus unfold more general  Thom-Boardman singularities  with controlled curvature  start by  observing  that  
 our image curve  $\eta_\circ(\mathbb R)\subset \mathbb R_+\times [-1,1]$, (which is 
  is only $C^1$-smooth),  is equal to the boundary of the  1-neighbourhood of the
   ray $ [1,\infty ) \subset \mathbb R\times [-1,]1]$.
  
  Then, to unfold  $\Sigma^{{1,...,1}}$,  of   depth $k$, where 
${{1,...,1}}=\underset {k}{\underbrace {{1,...,1}}}$, the natural model 
to use is  the  boundary of the 1-neighbourhood of the positive 
quadrant $\mathbb R^k_+\subset \mathbb R^k\times [-1,1]$, which  
has $curv\leq 1$ as well. But I haven't checked if this actually  works.\footnote {Beware of  non-coorientable folds, such as of the 
M\"obius  strip along he central line.}

 (b) It could be  interesting to quantify the  
   approximation procedure   
  of {\it smooth maps by immersion in Sobolev spaces} from   [GE 1971']  and also  a similar approximation in  [Be1991].

 (c) It is unclear how to "controllably unfold" in  $\mathbb R^{m+l}$ more general singularities  of smooth maps $X^m\to\mathbb R^m\subset \mathbb R^{m+l}$.

This leaves the following question open.

{\sf Do smooth immersions $f:X^m\to \mathbb R^{m+l}$ are regularly homotopic to immersions $f_\circ$, the curvatures of which 
are bounded up to a   multiplicative constant  by the minimal relative  curvatures of 
$\nabla^\perp$-trivial  immersions of   flat  tori  $\mathbb T^m\to \mathbb R^{m+l}$}.\vspace {1mm}

For instance it remains {\color {red!50!black}problematic} if 

{\it all $m$-manifolds  $X$ admit immersions 
$f:X\to B^{2m}(1)$ with curvatures $curv_F^\perp(X)\leq cost \sqrt m$, say for $const=100$.}

\section {Lower Bounds  on Curvature and upper Bounds on Expansion }

\subsection {Briefly on Scalar  Curvature:  $\bf \nexists PSC$,  $\bf Sc^\rtimes$, $\nexists \bf pss$, $\bf Ros.ind, SYS, etc$}

    All known    lower  bounds on the   curvatures of immersions $X\hookrightarrow Y$ (except for hypersurfaces in spheres\footnote{   If $curv^\perp (X^m\hookrightarrow S^{m+1})< 1$, then $X^m$ is homeomorphic  to $S^m$.
See  [Ge 2021] and section   3.7.3 in [Gr2021].}) depend on
   obstruction to positivity of the scalar curvature of Riemannian metrics  on $X$ and/or on submanifolds in $X$.

   Below is a (non-complete)  summary of what is known in this respect.\vspace{1mm}
      
  \textbf {$\bf Sc^\rtimes (X)$ and  $\bf \nexists PSC$.}
   Let   $X$  be a compact Riemannian manifold  (possibly) with a boundary, let $Sc(X,x)$  denote the scalar curvature at $x\in X$,
   that is  the sum of the values of the sectional curvatures $\kappa$ at the $m(m-1)$  (ordered) orthonormal  bivectors in  $T_x(X)$.
   
   For instance, $Sc(S^m(R))=m(m-1)/R^2$.

       Let $\lambda^{[\beta]}$ be  the lowest eigenvalue of the operator $$-\Delta_X+\beta Sc(X)$$ on $X$ with    the Dirichlet boundary condition. 
 
 Recall that
  $$\lambda^{[\beta]}(X)=\sup_\Theta\inf_{x\in X} \left ( Sc(X,x) -
  \beta^{-1}(\Delta \Theta(x)+\|\nabla \Theta\|^2)\right),$$ 
  where the supremum is taken  over all smooth functions $\Theta(x)$ on $X$.

  that $\lambda^{[1/4]}(X)$, denoted $$Sc^\rtimes (X),$$
    serves as a  worthwhile   substitute for $\inf_{x\in X}Sc(X,x)$ 
   (see [Gr2023]),  for compact as  well   as noncompact ones where in the latter case $Sc^\rtimes (X)$ is defined as the limit of $Sc^\rtimes (X_i))$
  for compact manifolds $X_1\subset X_2\subset ... \subset X$, which exhaust $X$.

{\it Examples.}  If  $Y$ has constant scalar curvature $\sigma$, then 
$$Sc^\rtimes(X)=4\lambda_1(X)+\sigma$$

For instance, the rectangular solids satisfy
$$Sc^{\rtimes}\left( \bigtimes_1^n[-a_i,b_i]\right)=4\sum_1^n \lambda_1[a_i,b_i]=
\sum_1^n \frac {4\pi^2}{(b_i-a_i)^2},$$
the unit hemispheres  satisfy:
$$Sc^\rtimes \left(S^n_+\right )=n(n-1)+4n=n(n+3)$$ 
 and 
 $$Sc^\rtimes (B^n)=4j_\nu^2,$$
 for the first zero of the Bessel function $J_\nu$,   $\nu=\frac {n}{2}-1$, where 
 $j_{-1/2}=\frac {\pi}{2}$, $j_0=2.4042...,$  $j_{1/2}=\pi$ and if $\nu>1/2$, then 

    $$\nu+\frac{a \nu^\frac{1}{3}}{2^\frac{1}{3}}<j_\nu<
\nu+\frac{a \nu^\frac{1}{3}}{2^\frac{1}{3}}+\frac{3}{20}\frac{2^\frac{2}{3}a^2}{\nu^\frac{1}{2}}$$
   where $a=\left(\frac {9\pi}{8} \right)^\frac{2}{3}(1+\varepsilon)\approx 2.32$  with
  $\varepsilon< 0.13\left(\frac {8}{2.847\pi}\right )^{2}<0.1$    [QW1999].

{\it Question.} Is there a  non-trivial bound 
$$\lambda_1(X_1\times [-r,r] )\geq  \lambda_1(B^{m+1} (1))+\varepsilon,$$
 where  $X=X^m$ is immersed to the unit ball 
$B^{m+1} (1)$  with $foc.rad(X)>r$, where  $ X_1\times [-r,r] $
is endowed  with the (flat) metric induced by the normal 
exponential  map  $ X_1\times [-r,r] \to B^{m+1}(1)$ and where 
$\varepsilon=\varepsilon (X)>0$ for non-spherical $X$, e.g.
$\varepsilon(X)\geq 1/10^{m+1}$ for  $X=_{homeo} \mathbb T^m$, $m\geq 2$?

    \textbf  {$\bf \nexists PSC$ and  Enlargeability.} 
  A smooth manifold  $X$ is $\exists PSC$ if admis a metric with $Sc>0$, otherwise it is called   $\nexists PSC$.
  
  For instance if a Riemannian manifold $X=(X,g)$ has $Sc^\rtimes >0$ 
  then it is $\exists PSC$. In fact, there exists a conformal metric 
  $Sc(\phi g)>0$  for some  function $\phi(x)>0 $ by he Kazdan-Warner theorem.

A  basic class of $\nexists PSC$ is constituted by {\it enlargeable} manifolds,\footnote {The implication $enlargeable  \implies \nexists PSC$ and related $\nexists PSC$ results and problems  are extensively discussed in 
    [Gr 2021], where the reader finds  further references. }
   where a closed  Riemannian  
  $m$-manifold $X=(X,g)$ is {\it enlargeable}  if it admits  a sequence  of orientable  covering $\tilde X_i\to X$  and   $\lambda_i$-Lipschitz maps 
  $f_i:(\tilde X_i, \tilde g_i) \to S^m(1)$,  such that 

$\bullet$ the maps$f_i$  are locally  constant at infinity and have {\it non-zero degrees}, 
  
 $\bullet$  $\lambda_i\to 0$ for $ i\to \infty$.}
 
 Clearly, enlargeability is a homotopy invariant of $X$; moreover,
if  $X_1\to X_2$ is a map with non zero degree and  $X_2$ is enlargeable, then $X_1$ is also enlargeable.
     
 {\it Examples.}
   Tori, as well as all manifolds  $X$   with $sect.curv^\perp(X)\leq 0$  are, obviously, enlargeable and (slightly less obviously)  {\it aspherical  locally homogeneous} Riemannian manifolds  are also enlargeable.

Aspherical $3$-manifolds are enlargeable  and  there is no example at the present moment  of a {\it non-enlargeable aspherical $m$-manifold}  for $m\geq 4.$

 {\it \textbf {Rosenberg Index.}}  This is an invariant  of closed  smooth manifolds,  which takes  values in some  (algebraic K-theory) group, where
  non-vanishing, $Ros.ind(X)\neq 0$  for spin manifolds $X$ is, essentially,  a shorthand for: 
{\sl
{ "$X$ is $\nexists PSC$, where this property is provable by a Dirac-theoretic argument"}.}
  
   {\it Examples.} The following  manifolds $X$, if spin, have    $Ros.ind(X)\neq 0$.
  
  $\bullet_1$    $4k$-Manifolds, where a certain Pontryagin number,   called   $\hat A$-{\it genus,} {\it doesn't vanish.}
     
   $\bullet_2$ {\it Hitchin’s spheres}:   manifolds   homeomorphic  (but not diffeomorphic)  to the spheres $S^m$,  $m = 8k + 1; 8k + 2$,  which don't bound spin manifolds; these do  exist for all $k = 1, 2, 3, ...$.     
     
    $\bullet_3$  Enlargeable manifolds and their  products  by 
    those in $\bullet_1\&\bullet_2$. 
    
    Moreover, if $\underline X^l$ is enlargeable, and $X^{l+m}$ 
    admits a smooth map $X^{l+m}\to \underline X^{l}$, such that  
     the pullback $Y=Y_x$ of a generic point $x\in X^l$ has $Ros.ind(Y)\neq 0$,  then  $Ros.ind(X^{l+m})\neq 0$  (see [WXY2021], [Ku 2021] and references therein.)
     
     \vspace {1mm}
    
  {\it $\widetilde  {spin}$-Remark.}   If  $Ros.ind(X)$ is {\it non-torsion}, then, in many (all?) cases
   (e.g. for enlargeable 
     manifolds and their products by those with $\hat A\neq 0$)
      the spin condition on $X$ in the poof of $\nexists PSC$  can be replaced by spin of the universal covering $\tilde X$  of $X$, where  the latter is satisfied, for instance, if $\pi_2(X)=0$. 
     
    \textbf  {$\bf SYS$-Manifolds,    $\bf SYS$-Enlargeability and 
     $\bf\nexists \bf pss$}.   SYS is  a condition on the integer homology of  $X$ introduced  by Schoen and Yau  who proved in [SY1979]   using {\it minimal hypersurfaces}   that $SYS\implies \nexists PSC$ (here $X$ doesn't need to be spin) 
 for $m=dim(X)\leq 7$, where this inequality is due to a
  possible  existence  of {\it \textbf perturbation \textbf stable \textbf singularities} of minimizing hypersurfaces.

{\color {magenta} Conjecturally}  $\nexists \bf pss$ holds for minimizing hypersurfaces   of all  dimensions (this means that the set of metric $g$ on $X$, such that all  $g$-volumes  minimizing hypersurfaces $\Sigma\subset X $ are smooth,
is $C^2$-dense in the space  of all Riemannian  metrics on $X.$   Also, $\bf \nexists pss$   is expected  of  {\it stable $\mu$-bubbles.} 

 This was confirmed for $m=8$ in [Sm1993] 
and in [CMS  2023] for $m\leq 10$ in the volume minimizing case.

Besides   a  {\it 2d-partial}   $\nexists \bf pss$ is presented  in  [SY2017],where it is used  for the proof of $\nexists PSC$ for SYS-manifolds of  all dimensions $m$.

Experience shows,    the types  of arguments used in  these papers for minimizing hypersurfaces 
equally apply to the  stable $\mu$-bubbles, but I checked this only for $dim(X)\leq 8$.

$SYS$-enlargeability generalizes enlargeability, for instance 

{\sl products of $SYS$-manifolds by enlargeable ones are $SYS$-enlargeable,}

{\sl circle bundles  over enlargeable manifolds {\sf(albeit not necessarily enlargeable)} are $SYS$-enlargeable}

and 

 {\it 2d-partial   $\nexists \bf pss$  suffices  for the implication 
$SYS$-enlargeable $\implies \nexists PSC$} [Gr2018]

\textbf { $\bf\nexists PSC\times Enlargeable$.} {\it If $X$  is    $\nexists PSC$ manifolds and $Z$ is a closed  enlargeable  manifold   of dimension  $m\neq 4$, then, {\color {magenta} granted $\nexists \bf pss$  for the stable $\mu$-bubbles}, 
the  product $X\times Z$ is $\nexists PSC$.}
 
Indeed,  $\bf\nexists  pss$  allows, for all $\varepsilon>0$,  a representation of the homology class $[X]\in H_m(X\times Z)$  by a submanifold $X_\varepsilon \subset  X\times Z$, which is (normally) frame bordant in $X\times Z$ o $X$\footnote {This, strictly speaking, makes sense only for orientable $X$, and it should be phrased more carefully if $X$ s non-orientable, where one must be aware that double covers of $\nexists PSC$ manifolds can be  $\exists PSC$. }   
and  such that 
$$Sc^\rtimes (X_\varepsilon)\geq Sc(X\times Z)-\varepsilon$$ 
(see section 2 in [Gr2023]).
  
  Then a  codimension $\geq 2$-surgery applied to  $X_\varepsilon $ brings you back to $X$, now endowed with a metric with $Sc(X)\geq Sc(X_\varepsilon-\varepsilon${\footnote {If $X$  is spin  this follows directly from  theorem1.8 in [St2002].}
  
  {\it Remark.} The simplest case of this, where $Z=S^1$ and where    
    $\bf\nexists  pss$ is needed only for minimal  hypersurfaces; thus the above is unconditionally true  for
     $dim(X\times Z)\leq10.$\footnote{In fact, the $\times S^1$-stability of   $\nexists PSC\times$ is formulated  as  conjecture  1.24  in [R 2006].
  
  But  since the above argument for $m=2,3,5,6$ is missing  in  this paper, I am worried    of myself making a silly mistake.}
    
   All   applications of
  $\bf\nexists pss$ in the present paper are derived from  the following. 
  
  \textbf {5.1.A. $\mathbb T^\rtimes$-Stabilized Band Inequality.} Let $X$  be a Riemannian manifold without boundary and let $f:X\to [-r, r]$  be a 1-{\it Lipschitz}, i.e. {\it distance non-increasing function.} Then, {\color {magenta}granted $\bf\nexists  pss$ for stable $\mu$-bubles of dimensions $\leq dim(X)-1$}, there exists a smooth  properly embedded hypersurface $\Sigma$, which {\it separates 
  $f^{-1}(-r)\subset X$ from} $f^{-1}(r)\subset X$\footnote {This means 
  that all  curve-segments in $X$ with the  ends in   $f^{-1}(-r)$ and  $f^{-1}(r)$ 
   intersect $\Sigma$ . (This condition is non-vacuous only if both sets $f^{-1}(-r)$ and  $f^{-1}(r)$ are non-empty.)} and such that
  {\color {blue}$$Sc^\rtimes(\Sigma)\geq Sc^\rtimes (X)-Sc^\rtimes ([-r, r])= Sc^\rtimes (X)-\pi^2/r^2.$$}
   (See  [Gr 2023] and references therein.)

   intersect $\Sigma$ . (This condition is non-vacuous only if both sets $f^{-1}(-r)$ and  $f^{-1}(r)$ are non-empty.)} and such that
  {\color {blue}$$Sc^\rtimes(\Sigma)\geq Sc^\rtimes (X)-Sc^\rtimes ([-r, r])= Sc^\rtimes (X)-\pi^2/r^2.$$}
   (See  [Gr 2023] and references therein.)





  \subsection{Norms on Curvature,   $m$-th   Scalar Curvature, Gauss Formula and Petrunin's Inequality}


  Besides   the normal curvature of an immersion $f:X\to Y$ at a point $x\in X $,  that is
   $$curv^\perp_f(X,x)=\sup_{ \tau, \nu}\|\rm II_\nu(\tau,\tau)\|, $$
   where II = II$_\nu(\tau_1,\tau_2)$  is the second fundamental form the supremum is taken over the unit tangent vectors $\tau\in T_x(X)$ and unit normal vectors $\nu\in  T^\perp_{f(x)}(X)=T_{f(x)}(Y)\ominus T_{x}(X)$,
  define the $l_2$-norm of the second fundamental form 
 ${\rm II= II}_f(X,x)$ 
 as follows.
    $$\|\rm II\|^2_{l_2}=
    \sum_{j=1,...k} \sum_{i_1, i_2=1,...m}\rm II_{\nu_j}(\tau_{i_1},\tau_{i_2}))^2,$$
  where  $\{\tau_i\}$, $i=1,...,m=dim(X),$   is a frame of  orthonormal vectors  in the tangent space $ T_x(X)\subset T_{f(x)}(Y) $ and $\{\nu_j\}$, $j=1,...,k=codim(X\hookrightarrow Y)$ is such a frame in the normal space   $T^\perp_{f(x)}(X)=T_{f(x)}(Y)\ominus T_{x}(X).$

  Observe that
   $$\|{\rm II}\|^2_{l_2}=
\sum^k_{j=1} \sum_{i_1,i_2}{\rm II}_{\nu_j}(\tau_{i_1},\tau_{i_2}))^2\leq
   k\sum_{i_1, i_2}({\rm II}_{\nu_j}(\tau_{i_1},\tau_{i_2}))^2\leq km\cdot curv_F^\perp(X)^2  $$
  and that
    $$\|{\rm II\|^2_{l_2}=
     \sum_{i_1, i_2}\sum_j(\rm II_{\nu_j}(\tau_{i_1},\tau_{i_2}))^2}\leq
   m^2\sum_j({\rm II_{\nu_j}(\tau_{i_1},\tau_{i_2}))^2}\leq  {m^2}curv_F^\perp(X)^2,$$
because
$$
curv_F^\perp(X)^2= \sup_{\nu, \tau} \|\rm II_\nu||(\tau,\tau)\|^2\geq
 \sup_{j,i_1,i_2}\| \rm II_{\nu_j}(\tau_{i_1},\tau_{i_2})\|^2,
$$
 $$ \sum_{i_1,i_2}({\rm II}_{\nu_j}(\tau_{i_1},\tau_{i_2}))^2=\sum_{l=1}^m\alpha^2_{j,l}$$
  for principal curvatures $\alpha_{j,l}$ of $X$ with respect to $\nu_j$ and
 $$\sum_j(\rm II_{\nu_j}(\tau_{i_1},\tau_{i_2}))^2=\sup_{\nu} (\rm II_{\nu}(\tau_{i_1},\tau_{i_2}))^2,$$
  since  $ \rm II_{\nu},(\tau_{1},\tau_{2})$ is
   a {\it linear function} in $\nu\in T_{f(x)}^\perp(X) $ for all  pairs 
   $\tau_{1},\tau_{2}\in T_x(X).$\footnote {Gilles Pisier explained to me that random (in a suitable sense)    II  have $\|{\rm II}\|^2\geq const\cdot km\sup_{\nu. \tau}\rm II_\nu(\tau,\tau) |$, that is  the inequality 
   ||I\hspace {-0.3mm}I$_f(X)$||$^2\leq km\cdot curv_F^\perp(X)^2$ is optimal  up to a multiplicative constant. }

  Next, following [Pet2023],   define ${\sf \Pi=\Pi}_f(X,x)$ as the average of 
 $\sum_j|\hspace{-0.2mm}|$||$_{\nu_j}(\tau,\tau\hspace{-0.2mm}|^2$ over the unit vectors $\tau\in  S^{m-1}_x\subset T_x(X)$.
  
  Clearly, 
  $${\sf \Pi}_f(X,x)\leq (curv^\perp_f(X,x))^2,$$ 
  where the equality 
  holds if and only if $||$II$||^2_{l_2}=||mean.curv^\perp(X,x)||$.
 (If $codim(X)=1$, this means that all principal curvatures  $X$ at $x$ are mutually equal.)

 Furthermore,  integration of $||$II$(\tau,\tau)||^2$, which is a 4th-degree polynomial in $\tau$, over the unit sphere $S^{m-1}_x$  
shows (as in  section 2 and in  [Pet2023]) that
  $$ \mbox {${\sf \Pi}={2\over m(m+2)} \big(||$II$||^2_{l_2}+{1\over 2}|\hspace{-0.2mm}|mean.curv^\perp|\hspace{-0.2mm}|^2\big)$}$$
   or $$\mbox {$||$II$||^2_{l_2}={m(m+2)\over 2}{\sf \Pi} -{1\over  2}|\hspace{-0.2mm}|mean.curv^\perp|\hspace{-0.2mm}|^2$}$$

  For instance, if  $m=dim(X)=1$ and $n=m+k=dim(Y)=3$ and the  principal curvatures of $X$ at $x$ are $ \alpha_1$ and   $\alpha_2$, then 
 $$curv^\perp(X,x)=\max( |\alpha_1|,|\alpha_2|),$$
  $$\mbox{ $||$II$||^2_{l_2}=\alpha_1^2+\alpha_2^2$},$$ 
   $$|\hspace{-0.2mm}|mean.curv^\perp|\hspace{-0.2mm}|=|\alpha_1+\alpha_2|$$
   and $${\sf \Pi}={1\over 4} ( \alpha_1^2+\alpha_2^2)+{1\over 8} (\alpha_1+\alpha_2)^2={3\over 8}( \alpha_1^2+\alpha_2^2)+\frac{1}{4} \alpha_1\alpha_2;   $$
   if
   $X=S^2\subset Y=\mathbb R^3$, where  $\alpha_1=\alpha_2=1$, this  makes 
${  \sf \Pi}=1$ as well.

   \vspace {1mm}
   
   Now let us turn to the curvature the ambient  Riemannian manifold  $Y$ of dimension $n\geq m $  and   define the function $Sc_{|m}(Y)$ on the tangent $m$-planes $T_y^m\subset T(Y) $ in  $Y$,   as the sum of the sectional curvatures $\kappa$ of $Y$ on the bivectors in  $T_y^m$  at $y$,   that is  the scalar curvature of submanifold $y\ni Y^m_y\subset Y$ tangent $T_y^m$, i.e. $T_y(Y_y^m)=T_y^m\subset T_y(Y)$ and having zero relative curvature in $Y$ at $y$, 
 $$ Sc_{|m}(Y, T_y^m)=Sc(T_y^m, y) =\sum_{i\neq j=1,...,n}\kappa(e_1\wedge \nu_j)$$  for a  frame of ortonormal  vectors $\tau_i\in  T^m_y$

In this terms,  the Gauss formula for the scalar curvature of $X\hookrightarrow Y$ reads:
$$\mbox {$Sc(X,x) = Sc_{|m}(Y,T_x(X))+ ||mean.curv^\perp(X,x)||^2-\|\rm I I\|^2_{l_2}$}, $$ 
 where by Petrunin's formula 
$$\mbox{$ ||mean.curv^\perp(X,x)||^2-||$ II$(X,x)||_{l_2}^2=||{3\over 2}mean.curv^\perp(X,x)||^2- {m(m+2)\over  2}{\sf \Pi}$}. $$

Hence, the inequality  $Sc_{|m}(Y)\geq \sigma_m$ implies that 
$$\mbox {$Sc(X) \geq \sigma_m-||$ II$(X,x)||^2$}.$$
Therefore 
$$Sc(X) \geq \sigma_m- km\cdot curv^\perp(X)^2\leqno [km]$$ 
for $k\leq m$
and 
$$Sc(X) \geq \sigma_m-m^2curv^\perp(X)^2.$$
for all $k$, 
where  Petrunin's formula yields the better inequality    
$$Sc(X) \geq \sigma_m-{m(m+2)\over 2}{\sf \Pi} \geq \sigma_m-{m(m+2)\over 2}   curv^\perp(X)^2.$$

It follows that 
{\it  if  the manifold $X$ is {\color {red!50!black}$\mathbf\nexists$}PSC, i.e. it admis no metric with $Sc>0$, then 

$$ curv^\perp(X)\geq \sqrt {\sf \Pi} \geq  \sqrt \frac {2\sigma_m} {m(m+2)}\mbox {   for all   $k$  and   }  n=m+k=dim (Y),\mbox { }  Y\hookleftarrow X,$$
and   
 $$ curv^\perp(X)\geq \sqrt \frac {\sigma_m}{km}\mbox {  for $k<m/2$}.$$}        \vspace {1mm}
 
                    \hspace {31mm} {\sc Examples and  Corollaries.}
                    \vspace {1mm}

 {\it Let $X$ be an $m$-dimensional {\color {red!50!black}$\mathbf\nexists$}PSC manifold, 
e.g. the $m$-torus $\mathbb T^m$ or Hitchin's exotic sphere} and let us indicate several examples of lower bounds on curvatures of immersions from $X$ to "small" manifolds $Y$.
\vspace {1mm}

 ({\Large$\ast$}$_{S^{n}}$)  {\sf Immersions $f$ from $X$ to the unit spheres $S^n$, $n= m+k,$  where $Sc_{|m}(S^n)=Sc(S^m)=m(m-1),$  satisfy
$$ curv_F^\perp(X)\geq \sqrt  \frac {m-1}{k}$$
and 
$$ curv_F^\perp(X)\geq\sqrt{\sf  \Pi}\geq  \sqrt  \frac {2 m-2}{m+2}\mbox  {  for all }k, $$}
where the latter inequality (1.2.(e) in [Pet 2023]) is sharp as it  is seen    with  
 $\mathbb T^m \subset S^{2N_{Hilb}}(1)$ in section 2.1, where   
 $\sqrt{\sf  \Pi}_{Hilb}= curv (\mathbb T^m)=\sqrt  \frac {2 m-2}{m+2}$.

(c) {\sf  Let $Y$ be the product, $Y=Y_1^{n_1}\times Y_2^{n_2}$, where the sectional 
curvatures of   the factors are bounded from  below by 
  $$sect.curv^\perp(Y_1)\geq 1\mbox {  and }sect.curv^\perp (Y_2)\geq-\kappa.$$
 For instance $Y_1=S^{n_1} $  and $Y_2= \mathbf H^{n_2}.$
 
 Let 
 $f$ be an immersion from  an $m$-dimensional  {\color {red!50!black}$\mathbf\nexists$}PSC manifold to $Y$ and let
 $$\frac { (m-n_2)(m-n_2-1}{ n_2(n_2-1)}\geq \kappa.$$}
  {\it Then the  curvature of $X$ is bounded from below by
$$curv_F^\perp(X)\geq  \sqrt \frac  { (m-n_2)(m-n_2-1)-\kappa n_2(n_2-1)}{m(n_1+n_2-m)}.$$}
 For instance, 
 
 {\it if $m=2n_2< 2n_1$ and $\kappa=-\frac {1}{2}$, then 
  $$curv_F^\perp(X)\geq  \sqrt \frac  { n_2-1}{4(n_1-n_2)}.$$}

  Indeed, the $m$-th scalar curvature of $Y$ is bounded from below by
 $$Sc(Y)\geq {m-n_2}{m-n_2-1}-\kappa n_2(n_2-1).$$

 \vspace{1mm}

 {\it  \textbf {(a)  On Extremality of Veronese.}}  (a) The smallest  curvature of  immersions of  closed  connected  non-spherical  manifolds to $S^N$ is, {\color {magenta}  probbaly,} that of the  {\it Veronese embeddings} (see 6.1)
$\mathbb RP^m\to S^{ \frac {m(m+3)}{2}-1}(1)$, where 
 $\sqrt{\sf  \Pi}_{Ver}= curv (\mathbb RP^m)=\sqrt  \frac { m-1}{m+1}$.

Let us show in this regard that  

{\it closed  connected non-spherical manifolds $X^m$   admit no immersions to $S^N(1)$ with curvatures} 
$$curv^\perp(X\to S^N(1)) <\sqrt{ 3\over 10} 
\left (=\sqrt {{1\over 3}-{1\over 30}}\right).$$

Indeed, let such an $X=X^m$   be  immersed to he unit   $N$-sphere, such that
$$\mbox {$curv^\perp(X^m\hookrightarrow S^N(1))\leq \delta$, $m\geq 2$,   $\delta<1/\sqrt 2.$}$$ 
Then the sectional curvature of the induced metric   in $X$ satisfies 
$$ 1-2\delta^2\leq  sect.curv^\perp(X)\leq 1+2\delta^2.$$

$(\bullet)$ If  $\delta = curv^\perp(X^m\hookrightarrow S^N)\leq 1\sqrt 2$, 
  then
  $$sect.curv^\perp(X^m)\geq 0\implies rank(H_\ast(X;\mathbb F))\leq    const_m \mbox { for all fields $\mathbb F$}.$$
 and the universal covering $\tilde X$  of $X$ satisfies 
 $$ diam(\tilde X)\leq \pi/\sqrt {1-2\delta^2}.$$

 $(\bullet\bullet)$ If $\delta\leq \sqrt {3\over 10}$, then 
  $$ \frac {1+2\delta^2}{1-2\delta^2}\leq 4$$
   and the universal covering of $X$  is diffeomorphic to $S^m$. 

$(\bullet\bullet\bullet)$  If $curv^\perp((X^m)\hookrightarrow S^N(1))\leq \delta$ then 
 $$  inj.rad(X)\geq \pi\sqrt \frac {1} {1+\delta^2},$$
and if $X$ is not simply connected, then the universal covering of $X$  satisfies
$$  diam (\tilde X)\geq {2\pi\over\sqrt {1+\delta^2}}.$$

Since the inequality
 $${\pi\over \sqrt {1-2\delta^2}}\geq {2\pi\over\sqrt {1+\delta^2}}$$
implies that $\delta\geq {1\over \sqrt 3}>\sqrt {3\over 10}$,  the proof follows.

{\it \textbf {(b) Rigidity Remark.} } Looking closer at the above argument reveals  that 
\vspace {1mm}

{\sf the inequality 
$$curv^\perp (X^2\to S^N(1))\leq {1\over \sqrt 3}$$ 
for   closed connected {\it non-simply connected}   surfaces $X^2$ implies that 
these are} {\it congruent to the Veronese surface}.\vspace {1mm}

{\it \textbf {(c)  Small Volume Remark.}} If 
$$\mbox {$curv^\perp(X^m\hookrightarrow S^N(1))\leq \delta<1/\sqrt 2.$, $m\geq 2$,}$$
then the   volume of $X$ is, on one hand,  bounded by 
$$vol(X)\leq {\left ({1\over \sqrt {1-2\delta^2}}\right )^m} vol(S^m(1))$$
and on the other  hand
$$vol(X)\geq {\left ({1\over \sqrt {1+2\delta^2}}\right )^m}vol(S^m(1)).$$
It follows  that  {the order of the fundamental group of $X$ is {\it  bounded  by the ratio of these two numbers. 
$$card(\pi_1(X)) \leq \left ({\sqrt {1+2\delta^2}\over \sqrt {1-2\delta^2}}\right )^m.$$}

{\it \textbf {(d) On Almost Flat $X\hookrightarrow S^N$}.} If 
$curv^\perp (X\hookrightarrow S^N(1))\leq \delta$ for a small $\delta $ - I roughly checked this for $ \delta\leq 1/4$ - the manifold $X$  

{\it lies 
$5\delta$-close to an equatorial $S^m\subset S^N$  and the normal projection $X^m\to S^m$ is a  diffeomorphism.}

 {\it \textbf {(e) On Non-spherical  $Y\hookleftarrow X$}.} The inequalities   ({\Large$\ast$}$_{S^{n}}$) implies similar inequalities 
for Riemannian manifolds   $Y^N$ with boundaries  (e.g.  Euclidean  and hyperbolic balls)    by means of suitable (e.g. projective  as in the proof of  Burago's inequality in 3.2.3 in [Gr2086])  diffeomorphisms $F$ from   $Y^N$ 
into $S^N(1)$, where  these $F$ controllably increase  curvatures of curves in  $Y^N$.

For instance

 {\sf let  $Y$ be a complete Riemannian $n$-manifold with the sectional curvature  bounded in the absolute value by 
$|sect.curv^\perp(Y)|\leq 1$ and with $inj.rad\geq 1$  and let 
 $f$ be an immersion from  an $m$-dimensional  {\color {red!50!black}$\mathbf\nexists$}PSC manifold to $Y$ such that the diameter of the $f$-image of $X$ in $Y$ is at most  $0.1$.}
 
{\it Then the curvature of $X$ in $Y$ is bounded from below by 
$$curv_F^\perp(X)\geq 0.1 \sqrt \frac {m}{n-m} -1.$$}

{\it Sketch of the Proof.} Conformally modify the Riemannian metric $g$ of $Y$ in the vicinity  of $f(X)$ 
with a conformal factor $\phi(x)=\psi(dist_g (x,x_0))$, where $\psi(d)$ is defined by the following condition.

{\sf if $Y$ is isometric to the hyperbolic $n$-space $\mathbf H^n$ with curvature $-1$, then the unit  ball 
 $B_{x_0}\subset \mathbf H^n$, when endowed with the new metric
  $\psi(dist_{hyp} (x,x_0))g_{hyp}$ becomes  isometric to the
  hemisphere $S^n_+\subset S^n$.}

Then one can show that the curvature of $X$ with respect to the metric $\phi g$ is not much great than that that with respect to $g$ and he proof follows.

{\it Remark.} The  constant $0.1$ is very crude;  the proof of  a
similar nearly  optimal inequality is presented in the next section
following Petrunin's   argument  from \S 4 in [Pet2023].
\subsubsection 
     {On  Meanings of $Sc_{|m}$} 
     
     The above 
     notwithstanding,      the  role of  $Sc_{|m}$  and especially  of the bound $Sc_{|m}\geq \sigma$   played
 in shaping   the  global (geo)metric  and/or topological properties of Riemannian $n$-manifolds $Y$ with $Sc_{|m}(Y)\geq \sigma $ 
for $n>m\geq 3$ essentially   remain 100\% problematic.

 One expects that positivity of $[Sc_{|m}](Y)$  for $m<n=dim(Y)$  has greater  
significance than positivity of Sc(Y)= $[Sc_{|m}](Y)$. Below is
 an, albeit  weak,   confirmation to this.\footnote{Much stronger results   for another "intermediate  scalar curvaure"  are obtained in [BHJ 2022].}   

  {\sf Let $Y$ be  a Riemannian manifold, 
the boundary $\partial Y$ of which is divided into two disjoint parts,  $\partial Y=\partial_-Y\sqcup \partial_+ Y$, where $\partial_\pm Y$ are unions of connected components of $\partial Y$.}

{\it Let  
 $$ dist(\partial_-Y, \partial_+ Y)=2r,$$
let the  sectional curvature of $Y$ be bounded from below,
   $$\kappa(Y)\geq\kappa_-$$
and let 
$$Sc_{|(n-1}\geq \sigma.$$ }

\hspace {-6mm} Then 

{\it $Y$ contains a smooth hypersurface $X\subset Y$, which separates 
$\partial_- Y$ from $\partial_+ Y$} ({\sf recall that $\partial Y= \partial_- Y \sqcup\partial_+ Y$}) 
{\it and such that the scalar curvature of the induced Riemannian metric in $X$ satisfies:
$$Sc(X)\geq \sigma- (n-1) \alpha_{\kappa_-}(r)^2,\leqno  {\color {blue}[\sigma |\alpha]}$$}
where $\alpha_{\kappa_-}(r)$ denotes the curvature of the circle of radius $r$ in the
 standart  surface
 with constant curvature $\kappa_-$, e.g. 

$\bullet $    $\alpha_{1}(r)= \frac {\cos r}{\sin r},$

$\bullet $   $\alpha_{0}(r)= \frac{1}{r},$

$\bullet $ $\alpha_{-1}(r)= \frac {e^r +e^{-r}}{e^r -e^{-r}}$.

{\it Proof.} Let $X_{[2r]}\subset Y$ be the $2r$-equidistance hypersurface to $\partial_-Y$  and 
$X_{[2r-r]}\subset Y$ be the $r$-equidistant to $X_{2r}$ on the side of $\partial_- Y$.
Then clearly

({\Large \color {blue}$\circ$}$_r$)  {\it the hypersurface $X_{[2r-r]}$ is $C^{1,1}$-smooth with the curvature, i.e. with  the norm of the second fundamental form, bounded by  $\alpha_{\kappa_-}(r)$.}

Hence, $X_{[2r-r]}$ can be approximated by $C^\infty$-smooth hypersurfaces $X_\varepsilon\subset Y$  with curvatures bounded by $\alpha_{\kappa_-}(r)=\varepsilon$  for all $\varepsilon>0$. QED.
\vspace {1mm} 

{\it Remark.}  If $n\leq 8$ (and $\bf \nexists pss$ is known for the $mu$-bubbles),   then  $\partial_-Y$  and 
$X_2\subset Y$ can be separated by a  smooth stable $\mu$-bubble $X\subset Y$ such that 
the scalar curvature of a warped product metric $g^\rtimes= g^\rtimes(x,t) =dx^2+\phi(x)^2dt^2$ on $X\times \mathbb T^1$
is bounded from below in terms of $\sigma=\inf_ySc(Y,y)$  and $r$ as follows (see section 3.7 in
[Gr2021]),
$$Sc(X)\geq \sigma-\frac {(n-1)\pi^2}{nr^2 }.$$

Although this is not formally stronger than  {\color {blue}[$\sigma|\alpha$]},  it is  by far more general and informative.  
\vspace {1mm}

{\it Questions.} (a) Does
({\Large \color {blue}$\circ$}$_r$) generalize to
submanifolds $X\subset Y$ of codimensions $k>1$, where $Y$ is, in some way,  "wide  in $k$-directions"?

For instance, let  $Y$ be a Riemannian manifold homeomorphic  to $X_0\times B^k(1)$, where $X_0$ is a closed manifold of dimension $n-k$, let the sectional curvature  of $Y$ be bounded 
by  $|\kappa(Y)|\leq 1$ and the injectivity radius by $inj.rad(Y)\geq 1$ (compare with [Gr2022]).\vspace {1mm}

{\sf What else need you  know about $Y$ to effectively evaluate the minimal  $\alpha$, such that  $Y$ contains a submanifold $X\subset Y$ {\it homologous to} 
$X_0=X_0\times \{0\}\subset X_0\times B^k(1)=X$,
such that the  curvature of $X$ in $Y$ is bounded by $\alpha$?

What is the best bound on $\alpha$ in a presence of a {\it proper} (boundary-to-boundary)  {\it $\lambda$-Lipschitz} map $X\to B^k(1)$?}\vspace {1mm}

The known (unless I am missing some) {\it quantitative transversality theorems} applied to maps  
$X\to B^k$ deliver submanifolds  $X\subset Y$ with $\alpha\leq const_n$, but we need $X$ with
  $\alpha\leq const_k$  for our purposes.

Alternatively, an inductive use of ({\Large \color {blue}$\circ$}$_r$) leads to a bound
with 
$$const  \sim 100^{k(1+ diam(Y))}$$
 but this is not satisfactory either.  \vspace {1mm}

(b) {\sf How much (if at all) do (essential)  global (geo)metric  and/or topological properties of Riemannian $n$-manifolds $Y$ with $Sc_{|m}(Y)\geq m(m-1) $ 
for $m\geq 3$   differ from those with $Sc(Y)\geq n(n-1)$?}

For instance, 
does the product $\mathbb T^{n-2} \times S^2$, $n\geq 4$,  admit a metric with $Sc_{|3}>0$?


\subsection { $Mean_{|m}(\partial Y)$-Curvature, $Sc^\rtimes$-Curvature  and Immersions to Riemannian Manifols with Boundaries} 

Let  $\Phi=\Phi(y) $ be a smooth function on a Riemannian manifold $Y$, let $\nabla\Phi$  be the   gradient field of $\phi$ and 
let $\Delta_{|m} \Phi(y)$ be the maximum  of the $V$-derivatives of the $m$-volumes on the tangent $m$-planes $\tau\subset T_y(Y).$

{\it Example.} Let $Y=\mathbb R^N=\{y_1,...,y_N\}$ and 
and $ \Phi(y)=-(y_1^2+....+y_{N-l}^2)$. If $l\geq m$,
then  $\Delta_{|m} \phi(y)=0$ and if $l< m$ then $\Delta_{|m} \Phi(y)=-2(m-l)$.
In particular if $m=N$,  then   $\Delta_{|m} \Phi(y)=\Delta \Phi(y) $
where  $\Delta \Phi$ is  the  ordinary Laplacian.  

******************************************

Given an  $m$-submanifold $X^m\subset Y$, decompose $V$ on $X$ into a tangent and normal fields  to $X$, 
$V_{|X}=\tau_V+ \nu_V$, where $\tau_V=grad_X\Phi(x$
and recall that  the $V$-derivative of the volume (element) on $X$ 
is
$$\Delta_X\Phi(x)+ \langle {\sf M}(x), \nabla\Phi(x)\rangle,$$
 where $\Delta_X$ is the Laplacian on $X$ with respect to the induced metric on $X$  and where ${\sf M}(x)$  is the mean curvature vector  of $X$, such that  $\mp$-direction of {\sf M}
is chosen  such  that the  volume of  $X$ {\it increases} under the $V$-flow.

It follows that 
 $$-\Delta_X\Phi\geq  -\Delta_{|m}\Phi+\langle {\sf M}, \nabla\Phi\rangle.$$
 
 This yields the following lower bound on   the first eigenvalue $\lambda^{[\beta]}$ of the operator $-\Delta_X+\beta Sc(X)$ (see section  5.1) on $X$ with the induced Riemannian metric and   the Dirichlet boundary condition. 
 $$\lambda^{[\beta]}\geq \inf_{x\in X} \left (M(x)^2-\|{\rm II}\|l_{l_2}^2 (x)
 - \beta^{-1}\left (\Delta_{|m}\Phi+\|\nabla \Phi(x)\|^2\right )+
 \langle {\sf M(x)}, \nabla\Phi(x)\rangle\right)$$ 
 $$\geq  \inf_{x\in X}  \left (-\|{\rm II}\|_{l_2}^2 (x)
 - \beta^{-1}\Delta_{|m}\Phi- \left (\beta^{-1}+\frac {1}{2}\right) \|\nabla \Phi(x)\|^2\right) $$
  for all  smooth  functions $\Phi$ on $Y$.

$\bf Mean_{|m}curv\geq m$.} Define the   mean curvature counterpart of $Sc_{|m+1}$ as 
 $$mean_{|m}curv^\perp(\partial(Y)=c_1+c_2+...+c_m,$$ 
where  $c_i$ are  the  first $m$  smallest  principal  curvatures of he boundary  of $Y$.

{\it Example.} If $Y=B^m(R)\times\mathbb R^l\subset \mathbb R^{k+l}$, then $mean_{|m}curv^\perp(\partial Y)=\max ((m-l)R^{-1} ,0)).$

(Our understanding of $mean_{|m}curv$ is as meager as that of $Sc_{|m}$.  For instance, 
{\sl does  the product   $ \mathbb T^{n-2}\times B^2$,  $n\geq 4$. admit an immersion to $\mathbb R^n$ with $mean_{|3}$-convex boundary, i.e. 
with $mean_{|3}curv^\perp(\partial T^{n-2}\times B^2)>0$?}

This may be related to the above  question about 
$Sc_{|m+1}$ by the doubling construction  as it is done for  convex ($m=1$) and mean convex ($m=n-1$)  hypersurfaces.)

Let $Y$ be a complete  connected  Riemannian manifold with non-negative sectional curvature and with non empty boundary,  and let
$$mean_{|m}curv\geq m.$$

Let $\Phi(y)=-\gamma(1- \frac {1}{2} dist(y, \partial Y))^2$
and observe that 
$-\Delta_{|m}\Phi\geq m$ and $\|\nabla\Phi\|\leq 1$.

Let  $codim(X)=dim(Y)-m\leq k$   and $curv^\perp(X)\leq \sqrt C$
then
$$\lambda^{[\beta]}\geq - Ckm+
  \beta^{-1}\gamma m- \gamma^2\left (\beta^{-1}+\frac {1}{2}\right) \mbox   { for all $\gamma>0$}.$$

The   relevant for us $\beta={1\over 4}$, where
$$\lambda^{[1/4]}=Sc^\rtimes (X)\geq \max_{\gamma>0} \left(4\gamma m- Ckm-
{9\over 2}\gamma^2\right)=\frac {16m}{9} -Ckm-\frac {8m}{9}$$
$$=\left (\frac {8m}{9} -Ck
\right) m.$$

For instance if  $Sc^\rtimes (X) )\leq 0$, then 
$$curv^\perp(X)\geq \sqrt{ \frac {8m}{9k}}.$$

\textbf {5.3.A. Corollary.} {\it  Let $f: \mathbb R^m\to B^{m+k}(1) \subset  \mathbb R^{m+k}$ 
be a smooth   distance increasing immersion.}  Then 
 $$curv^\perp_f( \mathbb R^m)\geq \sqrt{ \frac {8m}{9k}}.$$

\vspace{1mm}

{\it Remarks.} (a)  The above is a (minor) modification of (a part of)  Petrunin's  argument  from   [Pet2023] used for  the proof  of the inequality
  $$curv^\perp_f( \mathbb X\overset {f}\hookrightarrow B^N(1) )\geq \sqrt{ \frac {3m}{m+2}}.$$
for immersions of $\nexists PSC$ manifolds $X$.

(b) If the curvature of an immersion of a closed manifold to the unit ball  satisfies 
$$curv^\perp(X^m\hookrightarrow B^N(1))\leq \delta$$   
 then  $X^m$  is contained in the  $20 \delta$-neighbourhood 
of an equatorial $S_{eq}^m\subset S^N$ and if $\delta\leq 0.05$, then  the normal projection $X^m\to S_{eq}^m$ is a diffeomorphism.

(This can be derived  from 5.2(d) by radially projecting $X^m$ to the boundary sphere $S^{N-1}=\partial B^N(1)$.)

\textbf {5.3.B. Problem.}  Identify $m$-dimensional manifolds $X$, $m\geq 3$,  which admit  metrics with $Sc>0$, yet satisfy Petrunin's inequality 
$curv^\perp_f( X\overset {f}\hookrightarrow B^N(1) )\geq \sqrt{ \frac {3m}{m+2}}$
for all  $N$ and all immersions $f$.

  Moreover, classify $m$-manifolds $X$ which admit immersions $f:  X\hookrightarrow B^N$
with $curv^\perp_f( X\overset {f}\hookrightarrow B^N(1) )\leq (1-\delta)\sqrt{ \frac {3m}{m+2}}$ for a given  $\delta>0$.

\subsection {Low Bounds on Expansion of Equidimensional Immersions.}

Let $f:X\to Y$ be an  {\it expanding}.  i.e.  locally distance nondecreasing      map between  compact manifolds with boundaries

The simplest invariant which is  monotone increasing under  such an  $f$    is {\it the  inradius} of $X$,
 $$inrad(Y)\geq inrad(X),$$
where
$inrad(X)=\sup_{x\in X} dist (x,\partial X)$.

In fact, 
$$dist(x, \partial X)\leq  dist(f(x), \partial Y)\mbox 
{ for all  } x\in X,$$
since (some connected component  of)  the $f$ pullback of a curve from $f(x)$ to   $\partial Y$
in $Y$ connects $x$ with partial $X$.

This also shows  that if $$dist(x_1,x_2)\leq dist(x_1,\partial X),$$
 then 
$$dist (f(x_1), f(x_2))\geq dist (x_1), x_2),$$
that is $f$ is distance increasing on all balls in the interior of $X$.

 {\it Example of a Corollary.} Let $X\subset \mathbb R^n $ be the convex hull of two  balls  of  radii $r_1$ and $r_2\leq r_1$ in $\mathbb R^n$, such that  distance $d$  
between  their centers  is $\geq r_1$ and let $f:X\to \mathbb R^n$ be an expanding map. Then
$$diam(f(X))\geq 2r_1+\frac {r_1(d+r_2)}{r_1+d+r_2} +
\frac {r_2r_1}{r_1+d+r_2} =\frac {r_1(d+2r_2)}{r_1+d+r_2}$$

In fact, the image of $f$ contains the union of an $r_1$-ball  $B_1=B(r_1)\subset R^n$ and an $r$-ball $B_2=B(r)$ for 
 $r=\frac {r_1(d+2r_2)}{r_1+d+r_2}$, with the center  of $B_2$ in the boundary $\partial B_1$.

\textbf {5.4.A.} {\it Question.} Do  expanding self-maps $ X\to X$ 
of compact manifolds with boundaries send 
$\partial X\to \partial X$?

More interestingly,  let $X$ and$ Y$ be closed  connected  domains in the Euclidean space $\mathbb R^n$  and let $\lambda_1(X)$ and   $\lambda_1(Y)$  be the first eigenvalues of the Laplace operators in these domains with the Dirichlet boundary conditions.  

\textbf {5.4.B.  Large  RectangleTheorem.}  {\it If $X$  is a rectangular solid, $X=\bigtimes _{i=1}^n[0,d_i]$ and if 
$X$ admits an expanding map  $f:X\to Y$, then 
$$\lambda_1(X)\geq \lambda_1(Y).$$}

In fact, let $(g_f)$ be  the flat metric induced by$f$ in $X$ and observe that 

$( \bullet )$   $\lambda_1(X, g_f)\geq \lambda_1(Y).$    

$( \bullet  \bullet )$ there exists  a {\it distance decreasing} map $(X, g_f)\to X$ of  positive degree  (this is $f^{-1}$, of course).

This yields the proof by the following  theorem (see [Gr2023] and references herein).

\textbf {$Sc^\rtimes$-Extremality of 
 Rectangular Solids.}  {\sf Let   $Z$  be a compact orientable Riemannian manifold with a boundary and 
$$\Phi =\{\phi_1,...,\phi_n\}\to \bigtimes _{i=1}^n[0,d_i]= [0,d_1]\times...\times [0,d_n]$$
 be a continuous map, 
  where the functions $\phi_i:Z\to[0,d_i]$ are 1-Lipschitz} (distance non-increasing) and where $\Phi$  sends $\partial Z\to \partial \bigtimes _{i=1}^n[0,d_i]$.

{\it If the map $\Phi$ has non-zero degree, then
$$Sc^\rtimes (Z)\leq   Sc^\rtimes\left ( \bigtimes _{i=1}^n[0,d_i]\right).$$}

\textbf {5.4.C.} {\it Question.} Does \textbf {5.4.B}  holds true for all convex $X$?


\subsection {Bounds on Focal Radii of Immersions}

      {\it \textbf{ The focal radius }} of an immersed manifold    $X\overset {f}\hookrightarrow Y$, 
$$foc.rad(X)=foc.rad(X\hookrightarrow Y)= foc.rad_f(X)$$  
is the supremum of those $R$, for which 
 the differential of the  {\it normal exponential  map, denoted } 
$$\exp^\perp:T^\perp(X)\to Y,$$
is   {\it injective} along all normal  segments of length $< R$, where, in the case of a non-complete 
$Y$ or a presence of a boundary $\partial Y$, one has to say "{\it \color {red!20!black} defined  and injective}...".

 If $Y$ has constant sectional curvature, then the focal radii of submanifolds are  intimately  related to their curvatures  in $Y$.  
   
For instance, 
$$foc.rad(X\hookrightarrow \mathbb R^N)=\frac {1}{curv^\perp(X\hookrightarrow \mathbb R^N)}.$$
and 
$$foc.rad(X\hookrightarrow  B^N(1))=\min\left(\frac {1}{curv^\perp(X)}, dist(X, \partial Y)\right).$$

More generally,  if $curv^\perp(X)\leq \alpha$,  then $foc.rad(X)$  is bounded by the radii of circles in $S^2$ with  curvatures  $\alpha$ and if     $sect.curv^\perp(Y)\leq  \kappa$, then $foc.rad(X)$ is bounded
 from below by the radii of circles  in surfaces with constant curvature $\kappa$.

\textbf {   Codimension 1.}  If $X^m\hookrightarrow Y^{m+1}$ is a coorientable immersion,\footnote {Coorientability can be achieved by  taking  a double cover of $X$.}
 and $r< foc.rad(X)$
then the Riemannian metric $g_r$ induced in   
$ X\times [-r,r]$ by the exponential  map  satisfies   
$$\mbox {$Sc^\rtimes(X\times [-r,r])> Sc^\rtimes(Y)$ and $dist_{g_r}(X\times \{-r\},X\times \{r\})=2r$},$$
where the $\mathbb T^\rtimes$-stabilized {\it band 
 inequality} implies (see section 3.6 [Gr 2021] and references therein):
 $$r\leq  \frac { \pi}{ \sqrt {Sc^\rtimes(X\times [-r,r]) }}<   \frac { \pi}{ \sqrt {Sc^\rtimes(Y) }}
\leqno   
{\left[\frac { \pi}{ \sqrt {Sc^\rtimes}}\right]},$$
provided one of the following two conditions is satisfied.

(i) $X$  has  {\it non-zero Rosenberg index}\footnote{See section 5.1 for examples of manifolds  $X$ with $Ros.ind(X)\neq 0$.} 

(ii) No smooth hypersurface $X'\subset X\times [-1,1]$, which is homologous to $X\times \{0\}$ admits a metric with $Sc>0 $ 
 \footnote{If $m=dim(X)\neq 4$ this is equivalent to $\nexists PSC$=property of $X$.} and
  $m=dim(X)\leq 7$ (see section 5 in [Gr2021] and section 2 in [Gr2022]).\footnote{If $m>7$ one needs ,
    $\nexists \bf pss$  for the  stable bubbles in $(m+1)$-dimensional manifold, which, most likely, follows  for $m\leq 9$  by the argument from [CMS2023].}\vspace {1mm}

{ \it Comparison with the  Gaussian Curvature  Inequalities in Sections 5.2, 5.3}. 
If $Y=S^{m+1}(1)$, then the Gauss formula implies that 
$curv^\perp(X)\geq \sqrt {m-1}$ which shows that 
$$foc.rad(X)<{1\over  \sqrt {m-1}},$$
 while the above inequality  
$$foc.rad(X)<{\pi\over  \sqrt {m(m-1)}}$$
serves better for $m\geq 10\approx  \pi^2.$

Apparently,  there must be an inequality better than both of the two.

\vspace {1mm}

\textbf  {Expanding Quantification of $\left[\frac { \pi}{ \sqrt {Sc^\rtimes}}\right]$.} 
 Let $X$ be a compact orientable  $m$-dimensional manifold with a boundary,  let $Y$ be a Riemannian manifold of dimension $ m+1$ and let $X\to Y$ be a  smooth  immersion.

Let $\phi: X\to[-1,1]^m$ be a continuous map, which sends the boundary of $X$ to the boundary of the cube  $[-1,1]^m$ and let   $deg(f)\neq 0$.

Let  $d_i\geq r=2r$, $ i=1,...,m$ be 
 the distances  between the pullbacks of the opposite faces of the cube with respect the Riemannian metric  induced in $X$ by 
  the  immersion  $X\to Y$.

 Then, under the above assumptions  (i) or (ii), e.g. if $X$ is spin or if $dim(X)\leq 7$,  
$$ \frac {\pi^2}{4r^2} + \sum_1^m \frac {\pi^2}{(d_i-2r)^2} \geq \frac {1}{4} Sc^\rtimes (Y),$$
that is
$$r\leq \left( {Sc^\rtimes(Y)\over \pi^2}-\sum_{i=1}^m \frac {4}{(d_i-2r)^2}\right )^{-{1\over 2}}.\leqno {[d_i-2r]}$$

\textbf {Codimension 2.}  Let  be a  coorientable  immersion
 with $foc.rad (X)= r$. Then the exponential map 
 from the the normal ${r\over 2}$ -sphere bundle $S$ to $Y$  is 
an immersion with $foc.rad={r\over 2}$.

  \textbf {Case 1}: {\sf Euler$^\perp$} = 0. If $X$  is  $\nexists PCS$ manifold and  
  $dim(Y)=m+2=
dim(X)+2$, and the normal Euler class of the immersioN vanishes, 
(e.g.  $X\hookrightarrow Y$ 
 is an embedding homologous to $0$),  then 
$S\to X$ is a trivial circle bundle, $S= X\times S^1$.
Hence,  $S$ is $\nexists PCS$ with a few (some still unsettled)  exceptions 
(see section 5.1)  and 
$[\frac { \pi}{ \sqrt{Sc^\rtimes}}]$  applied to $S$ 
then shows that 
$$r\leq  \frac {2 \pi}{ \sqrt {Sc^\rtimes(X\times [-r,r])} } <\frac {2 \pi}{ \sqrt {Sc^\rtimes(Y)} }
\leqno   
{[\frac { 2\pi}{  \sqrt{Sc^\rtimes}}]}.$$
   and  all  smooth hypersurfaces $X''\subset X\times \mathbb T^1 \times [-1,1$, which are  homologous to $X\times  \mathbb T^1\times \{0\}$  satisfy either the above (i) or (ii).

  \textbf {Case 2}: {\sf Euler$^\perp\neq  0$.} Probably, if the Rosenberg index   is non-torsion,  then the same is true for  (non-trivial) circle bundles  $S\to X$ in most (all?)
  cases and  $[\frac { 2\pi}{  \sqrt{Sc^\rtimes}}]$ holds true in such a case  for spin manifolds.
  
 Also, if $X$ is enlargeable, then, granted $\nexists \bf pss$  for $m+2=dim (Y)$ for $Y\hookleftarrow X$, (unconditionally for $m= dim(X)\leq 6$), then 5.1.A  applies  and   shows that 
 the relative homology class of the fiber of the normal  $r$-disc bundle $D(r)\to X$ is realizable  by a surface $\Delta\subset D$
 with 
 $$Sc^\rtimes(\Delta)  \geq Sc^\rtimes(Y)$$ 
  and,  according to the $\mathbb T^\rtimes$-stable    Bonnet-Myers diameter inequality
 
   2.8(b) in [Gr2021],  
  $$r\leq \frac {2 \pi}{ \sqrt {Sc^\rtimes(Y)} }$$ 
 for enlargeable $X$ as well.

 (In fact, one needs here  only   a  {\it 2d-partial}   $\nexists \bf pss$ for the stable $\mu$-bubbles, while {\it 2d-partial}   $\nexists \bf pss$ for minimal hypersurfaces  yields the bound
  $r\leq 8 \pi/\sqrt{Sc^\rtimes(Y)} $.{\large )}.

\textbf {Codimension $\geq$  3.} Given a closed
orientable manifold $\Sigma$ of dimension $l\geq 2$, define {\it the  minimal} parametric (cospherical)  {\it area} 
$\bf PAR(\Sigma)$  
as

{\sl the infimum of the  numbers $A$}, such that there exist 

a compact smooth orientable $M$-dimensional manifold $P$ (possibly) with a boundary and a smooth     $P$-family $\Phi_p$  of smooth  maps from $\Sigma$ to the  unit  $l+M$-sphere, that is a smooth map
$$\Phi :\Sigma\to P\to S^{l+M}(1)$$ 
such that the maps  
$$\Phi_p:\Sigma \to S^{l+M}(1)$$ 
 are area {\it $A$-contracting} i.e.
 the areas of all smooth surfaces    $S\subset \Sigma$ 
$$area(F_p(S))\leq A^{-1} area(S)\mbox  {  for all } p\in P$$
and such that the map $\Phi$ is constant on  $\sigma\times \partial P$  and {\it the degree of $\Phi$ is  non-zero}.

.

One knows (see [Gr 20//] and references therein) that if  $\Sigma$ is spin,{\footnote{ The spin condition can be relaxed  
to that for the universal covering of $\tilde \Sigma$, but 
dropping spin all together  remains problematic.}  then
$$PAR(\Sigma)\cdot Sc^\rtimes(\Sigma) \leq PAR(S^l(1))
 Sc(S^{l})=4\pi l(l-1). \leqno {[Par<...]}.$$

Now, let $X$ be closed {\it enlargeable} $m$-manifold, $Y$ be a  Riemannian $(m+k)$-manifold,  let  $f:X\to Y$ be a smooth immersion, such that 
$$foc.rad_f(X)>r$$map from $Y$.
and let  $B=B(r)\to X$ be the  normal $r$-ball bundle endowed  with he Riemannian metric induced by the normal exponential 
map  $\exp^\perp: B\to Y$.

Then, granted 
  $\nexists \bf pss$ for $\mu$-bubbles of dimensions $\leq m+k-1$,\footnote{If 
   $k=3$, then  2d-partial   $\nexists \bf pss$ suffices.} there exits 
   a $(k-1)$-dimensional  submanifold $\Sigma=\Sigma^{k-1}$  in the 1/3-annulus in the  normal $r$-ball bundle   $B=B(r)\to X$
  $$\Sigma\subset B(2r/3)\setminus B(r/3)$$ with 
 $$Sc^\rtimes (\Sigma)\geq \sigma_r = 
 Sc^\rtimes (Y)-Sc^\rtimes [0,r/3]=Sc^\rtimes (Y)-36\pi^2/r^2.$$ 
and such that the homology class  $[\Sigma]\in H_{k-1}( B(2r/3)\setminus B(r/3)$  is equal to a {\it  non-zero multiple} of the class
of the fiber 
$S_x^{k-1} =\partial B_x(r/3)\subset B(2r/3)\setminus B(r/3)$.
(see [Gr2022] and references therein.

To simplify, let the normal bundle of $X$ in $Y$ be trivial\footnote {The vanishing of he Euler  class suffices  for our argument.}
let $P=X\times [r/3,2r/3]\hookrightarrow B(2r/3)\setminus B(r/3)$ where the embedding is naturally naturally associated with a section
of the $r$-sphere bundle $S=S(r)=\partial B(r)\to X$.

 Let $$B_p^\bullet(\rho)\subset B, \mbox { } p\in P, \rho\leq r/3,$$    be the one point 
 compactifications of the (open)  $\rho$-balls in $B$. 
 
 Let
 $$inj.rad_{\exp^\perp p} (Y)\geq r/3\footnote{ The  inequality $inj.rad_y(Y)\geq R$ signifies that $dist (y,\partial Y)\geq R$ and  that the exponential map $\exp_y$ in $Y$ at $y$ smoothly {\it embed} the tangent $R$-ball from $T_y(Y)$   to $Y$}$$
 choose a normal frame over $X$
and  (radially  diffeomorphically)  map  the  balls 
  $B_{\exp^\perp p}^{m+k} (\rho)\subset Y$ to $B^{m+k} (\rho)\subset \mathbb R^{m+k}$  by means of the inverse exponential   maps in $Y$  
 at the points $y=\exp^\perp p\in Y$.
 
 Finally, scale   $\bar B^{m+k} (\rho)$ to the $\pi$-ball $\bar B^{m+k} (\pi)$ and radially map it to 
  the unit sphere $S^{m+k}(1).$

 Thus we obtain  a map $F:\Sigma\times P\to S^{m+k}(1)$
 such  that
 
 if $\rho< r/3$ then $deg(F)=1$, 
  
 if $sect.curv^\perp(Y)\leq 1/\rho$, then $F$ 
  is $A$-contracting for 
 $A\geq \rho^2/\pi^2$. 
 
 Therefore, 
 
 {\it if $inj.rad_X (Y)\geq r/3$ and  
 $sect.curv^\perp(Y)\leq 1/3r$}, where $r=foc.rad(X)$,
 then
$$ (r^2/\pi^2)\cdot  (Sc^\rtimes (Y) -36\pi^2/r^2) \leq 4\pi(k-1)(k-2),$$
that is,
$$foc.rad (X)\leq \pi \sqrt {4\pi(k-1)(k-2) +36
\over Sc^\rtimes (Y)}<$$}

{\it Remarks.} (a) A similar result holds for immersions of product manifolds $X=X_0^l\times Z^m\to Y^{l+m+k}$, where  $Z=Z^m$ is enlargeable,
$$foc.rad (X)\leq \pi \sqrt {4\pi(k+l-1)(k+l-2) +36\over Sc^\rtimes (Y)}, 12 (k-1)$$
provided   $Euler (T^\perp(X))=0$ and granted  
  $\nexists \bf pss$ for $\mu$-bubbles of dimensions $\leq m+k+l-1$ \footnote{If $Z=\mathbb T^m$, then  $\nexists \bf pss$
is needed  for minimizing hypersurfaces of dimensions  
 $\leq m+k+l-1$, while  $\nexists \bf pss$ for the stable  $\mu$-bubbles  is needed only for  dimensions $\leq k-1$. And if $k=3$,  then    2d-partial   $\nexists \bf pss$ suffices.} 
 
(See [ Gr///] for  sharper and more general inequalities of this type.)

\subsection{Conjectures and Problems}
 
\textbf { Codimension $k$ Conjecture}. 
 {\it The inequality  $curv^\perp( X^m\hookrightarrow B^{n}(1)) \geq\frac{2j_\nu}{(n-m)\pi }-1$  holds for all compact 
 enlargeable $m$-manifolds,all  $n>m$ and    the first zero of the Bessel function $J_\nu$,   $\nu=\frac {n}{2}-1$.}

 \vspace {1mm}


 \textbf  {(Overoptimistic?) Conjecture.} {\sf If the cohomology  of a closed $m$-manifold $X^m$ with coefficients in some field  $K$
 contains $l$ elements  with non-zero  product,
 $$h_1,\smile...\smile h_i\smile ...\smile h_l\neq 0,\mbox { }  h_i\in H^\ast(X;K),$$
(e.g.  $X^m=S^{m_1}\times ...\times S^{m_l}$, $m_1+...+m_l=m$),} 

{\sl then the  curvatures of immersion $f:X^m\to B^{m+k}(1)$  bounded from below as follows,
$$curv_F^\perp(X) \geq 0.1 \frac {l^2}{mk}?$$}

\textbf {Clifford Tori Extremality Problem.} {\sl Does the $m$-torus admit an immersion to the unit $2m$-ball 
with curvature $<\sqrt m$?}
 
For all we know, all flat $m$-tori admit smooth isometric immersions to 
$B^{2m}(1) $ with curvatures < 10.

\textbf {$\mathbf {m^\beta}$-Problem,} {\sl What is the minimal $\beta$, such that   the tori  of all dimensions $m$  admit  immersion to the 
unit $(m+1)$-balls, 
$$f: \mathbb T^m\hookrightarrow B^{m+1}(1),$$
with curvatures $curv_F^\perp(\mathbb T^m)\leq 100 m^\beta?$ }
(We know that $\beta \leq \frac {3}{2}$.)

\textbf {Simply Connected Codim 1 Curvature Problem.} {\sl Do  all compact smoothly imbedded {\sf simply connected} hypersurfaces
 $X^m\subset \mathbb R^{m+1}$, e.g. products of spheres of dimensions $\geq 2$,  admit  immersion to the unit ball,  
$$f:X^m\hookrightarrow B^{m+1}(1)$$
with curvature $curv_F^\perp(X)\leq 100?$}


\textbf   A.  {\sf How large can  be  the ratio  
$$min.curv^\perp(X^m\hookrightarrow   B^{M}(1))/min.curv^\perp(X^m\hookrightarrow   B^{M+1}(1))$$
 provided $X$ immerses to  the Euclidean space $ \mathbb R^M$, e.g. for $M\geq 2m-1$?

Is this ratio  bounded by a universal constant, say by $const\leq100$?
 
\textbf   B. What is the {\it homological Morse spectrum} of the function 
$\mathcal M :f\mapsto  curv_F^\perp(X)$ on the space of immersions $f:X\to Y$?}

(An  $r\in \mathbb R_+$ is in the Morse spectrum of a function $\mathcal M:\mathcal F\to \mathbb R_+$  if there exists a homology class $h\in H_\ast(\mathcal F;A)$ with some coefficient group  $A$, such that  the $r$-sublevel of $\mathcal M$ contains $h$, i.e.  
$h$ is contained in the image of the inclusion homomorphism 
$$H_\ast(\mathcal M^{-1}[0,r];A)\to  H_\ast(X;A)$$
while  the lower sublevels $\mathcal M^{-1}[0,p]\subset \mathcal M^{-1}[0,r]$, $p<r$,  don't  contain $h$. See [Gr1988], [Gr2017] for more about it.)

\textbf   C. {What are bounds on the averages of powers of the curvatures of immersions,  
$$\frac {1}{vol_X}\int_X(curv_x(X\hookrightarrow Y))^pdx\mbox {  
for } p\geq 1?$$}
an 
 (See [Pet 2023] for such an inequality with the mean curvature 
 and consult [LB2021] for the  bounds on the Yamabe invariant  for 4-manifolds, which may (?)  apply here.

\textbf { Product of  Balls Problem}. Given positive numbers $r_i$,  $R_i$  and positive  integers $m_i$,
 $n_i$, $i=1,...k$,
such that   $\sum_im_i=\sum_in_i$, evaluate, let it be only roughly,  the maximal $\lambda>0$ , such that 
the product of  $m_i$-dimensional $r_i$-balls $B^{m_i}(r_i)\mathbb R^{m_i}$  admit a $\lambda$-expanding map 
to the product of $n_i$-dimensional $R_i$-balls,
$$\bigtimes _{i=1}^kB^{m_i}(r_i)\to \bigtimes _{i=1}^kB^{n_i}(R_i).$$

\textbf {Cube Extremality Problem.} {\sl Does, the unit $n$-cube $[-1,1]^n$ admits an expanding map to the $n$-ball  of radius $<\sqrt n$?}
\vspace {1mm}

\textbf {Expansion on  Mesoscale.} Find unified generalizations of the above results to classes of continuous maps $f:X\to Y$ stable under  $C^0$-perturbations.

\textbf B. {\it Example} 1. Given a function $\delta(d)$, study  continuous  maps
$f:X\to  Y$, such that 
$$dist_Y(f(x_1)f(x_2))\geq \delta (dist_X(x_1,_2)).$$  
(A possible $\delta$  may be  supported in the segment  $ [c, 100c]\ni d$, where 
  $\delta(d)\geq d $  for $c<d <100d$ and where eventually
    $c\to 0$.

\textbf C.  {\it Example} 2. Study embeddings $\phi:X^m\hookrightarrow V^n$ where $V\supset X$ retracts to $X$ 
and where $dist(X,\partial V)\geq r$;  then study composed maps 
$f= \psi\circ \phi:X\to Y^n$,
$$X^m\overset {\phi }\hookrightarrow  V \overset {\psi} \to Y^n.$$
where $\psi$   is an expanding map.

\textbf D.  {\it Example} 3.  Specialize the above to 
(piecewise linear) maps including  non-locally trivial $p.l.$ immersions and to more general piecewise smooth maps.

\vspace {1mm}
 \vspace{1mm}

 \section{Miscellaneous}
 

\subsection {Veronese Maps.}

Besides   invariant  tori, there are other submanifolds in the unit sphere   
$S^{N-1}$,  which have  small curvatures  and which are transitively acted upon by  subgroups in
the orthogonal   group $ O(N)$.
\vspace {1mm} 

{\it The  generalized  Veronese maps}   are a {\it minimal equivariant isometric}  immersions of spheres to spheres, with respect to certain homomorphisms ( representations)  between the orthogonal
groups $ O(m+1)\to O(m+1)$,
 $$ver=ver_s=ver_s^m:
 S^m(R_s)\to S^m=S^{m_s}= S^{m_s}(1),  $$  
 where  
 $$m_s= (2s+m-1)\frac{s+m-2)!}{s!(m-1!}<2^{s+m}\mbox { and } R_s=R_s(m)=\sqrt \frac {s(s+m-1)}{m},$$
for example, 
$$\mbox {$m_2= \frac{m(m+3)}{2}-1$,    $R_2(m)= \sqrt\frac {2(m+1)}{m}$   and  $R_2(1)=2$,}$$\vspace {2mm}
 (see [DW1971]
If $s=2$ these, called  {\it classical Veronese maps}, are defined  by taking squares of linear
functions (forms) $l=l(x)= \sum_i l_ix_i$  om $\mathbb R^{m+1}$,
$$Ver: \mathbb R^{m+1}\to \mathbb R^{M_m},\mbox { } M_m=\frac {(m+1)(m+2)}{2},$$
where tis   $\mathbb R^{M_m}$ is represented by the space $ \mathcal Q=\mathcal Q(\mathbb R^{m+1})$  of quadratic functions (forms) om $\mathbb  R^{m+1}$,
$$Q =\sum^{m+1,m+1}_{i=1,j=1}q_{ij}x_i x_j.$$

  The Veronese map, which is (obviously)    equivariant  for
the natural action of the orthogonal group  group $O(n+1)$ on  $ \mathcal Q$, where, observe, 
 this action fixes the line 
$\mathcal Q_\circ$ spanned  by the form $Q_\circ=\sum_ix^2$ as well as the complementary subspace $\mathcal Q_\diamond$ of the {\it traceless forms $Q$}, where the action of $O(n+1)$ is irreducible and, thus, it has  a {\it unique, up to scaling} 
Euclidean/Hilbertian structure. 

Then  the normal projection\footnote {The splitting $\mathcal Q=\mathcal Q_\circ \oplus \mathcal Q_\diamond$ is necessarily normal for all $O(m+1)$-invariant  Euclidean  metrics in $\mathcal Q$.}
defines an equivariant map to the  sphere in $\mathcal Q_\diamond$  
$$ver: S^m\to S^{M_m-2}(r)\subset \mathcal Q_\diamond, $$
where the radius of this sphere, a priori, depends on the normalization of the $O(m+1)$-invariant  metric in $ \mathcal Q_{\diamond}$. 

Since we want the map to be   isometric,  we  either take $r=\frac {1} {R_2(m)}=\sqrt\frac {m}{2(m+1)}$  and keep $S^m=S^m(1)$ or 
 if we let $r=1$ and  $S^m=S^m(R_2(m))$ for 
$R_2(m)=\sqrt\frac{2(m+1)}{m}$.

Also observe that the
Veronese maps,   which are not  embeddings themselves,  factor via embeddings of  projective spaces 
to spheres 
$$S^m\to \mathbb RP^m \subset S^{M_m-2}\subset
 \mathbb R^{M_m-1}=\mathcal Q_{\diamond},\mbox { } M_m=\frac {(m+1)(m+2)}{2}.$$

{\it \textbf {Curvature of Veronese.} } Let is show that

$$curv^\perp_{ver} \left(S^m(R_2(m))\hookrightarrow S^{M_m-2}(1)\right)   =\sqrt { \frac {R_2(1)} {R_2(m)}-1}=\sqrt\frac  {m-1}{m+1}.$$

Indeed, the Veronese map sends equatorial circles from $S^m(R_2(m))$ to planar circles of radii 
$R_2(m)/R_2(1)$, the curvatures of which in the ball $B^{M_m-1}$ is 
$R_2(1)/R_2(m)=2\sqrt\frac {m}{m+1}$  and the curvatures of these  in the sphere,  
$$curv^\perp(S^1\subset S^{M_m-2}(1))=\sqrt {curv (S^1\subset B^{M_m-1}(1))^2-1}=\sqrt {\frac {4m}{m+1}-1}=\sqrt \frac {3m-1}{m+1}$$
is equal to the curvature of  the Veronese $S^m(R_2(m))\hookrightarrow S^{M_m-2}(1)$ itself

$\sqrt{R_2(1)/R_2(m)}=\sqrt\frac {2m}{m+1}$, and the curvatures of these  in the sphere,  
$$curv^\perp(S^1\subset S^{M_m-2}(1))=\sqrt {curv (S^1\subset B^{M_m-1}(1))^2-1},$$
is equal to the curvature of  the Veronese $S^m(R_2(m))\hookrightarrow S^{M_m-2}(1)$itself.
QED.

\vspace {1mm}

It may be hard to prove (conjecture in section 1)   that  {\sf  Veronese manifolds  have  the smallest possible curvatures among  non-spherical $m$-manifold in the unit ball:}
{\sl if a smooth compact $m$-manifold $X$ admits a smooth  immersion to the unit ball  $B^N=B^N(1)$ 
 with curvature $curv^\perp(X\hookrightarrow B^N)<\sqrt \frac {2m}{m+1}$, then $X$ is diffeomorphic to $S^m$.}

It is  more realistic   to show that the Veronese  have smallest curvatures among  submanifolds $X\subset B^N$ {\it invariant under  subgroups in $ O(N)$, which  transitively act on $X$. }

\vspace{1mm} 

{\it Remark.} {\sf  Manifolds $X^m$ immersed to $S^{m+1}$ with curvatures $<1$ are diffeomorphic to $S^n$}, see 5.5,
   but, apart from   Veronese's, we {\color {red!60!black} can't rule out } such $X$  
 in  $S^N$ for $N\geq m+2$ \footnote { Hermitian Veronese maps from the complex projective spaces $\mathbb CP^m$ to the spaces $\mathcal H_n$ of  Hermitian forms on $\mathbb C^{m+1}$  are among the prime suspects in this regard.
 } and,   even less so,  non-spherical $X$ immersible with curvatures $<\sqrt 2$ to $B^{N}(1)$, even for $N=m+1$.

It seems hard to decide this way or another, but it may  be realistic to try to prove  {\it sphericity  of simply connected} manifolds immersed with curvatures $<1$  to $S^N(1)$   for all $N$.

\vspace {1mm}

The curvatures of Veronese maps can be  also evaluated with the {\it Gauss formula,} (teorema egregium),
  which also gives the following formula for curvatures of all $ver_s$:

$m=2$  $1-2c^2=1/3$, $2c^2=2/3$ $c\sqrt{1/3}$

$C=\sqrt{1+1/3}=2/\sqrt 3$

{\it \textbf {From Veronese to Tori.}}  The restriction of the  map 
$ver_s:S^{2m-1}(R_s)\to S^{N_s}$ to the Clifford torus $\mathbb T^m\subset S^{2m-1}(R_s)$  obviously satisfies 
$$curv^\perp_{ver_s}(\mathbb T^m)\leq A_{2m-1,s}+\frac {\sqrt m}{R_s}= 
\sqrt {3-\frac {5}{2}m+\varepsilon(m, s)}$$
 for 
 $$\varepsilon(m, s)= \frac {2}{4m^2}-\frac {4m-2}{s(s+2m-2)}+\frac{5(2m-1)}{2ms(s+2m-2)}-
  \frac {2m-1}{(ms(s+2m-2))^2}.$$

 This, for $s>>m^2$, makes $\varepsilon(m, s)=O\frac{1}{m^2}$

 Since $N_s< 2^{s+2m}$, 
 
 {\it starting from $N=2^{10m^3}$
 $$curv^\perp_{ver_s}(\mathbb T^m)<\sqrt{ 3- \frac {5}{2}m}.$$}
 where it should  be noted that

 {\it  the Veronese maps restricted to the Clifford tori are $\mathbb T^m$-equivariant}\vspace {1mm}

and that

{\sf this bound is {\it  weaker than  the optimal one}  $ \frac{||y||^2_{l_4}}{||y||^2}\geq\sqrt {3-\frac {3}{m+2}}+\varepsilon$ from 
 the previous section}.\vspace {1mm}

{\it Remarks}. (a) It is not hard to go to the (ultra)limit for $s\to \infty$  and thus obtain an 

{\sf equivariant  isometric immersion $ver_\infty$ of the 
Euclidean space $\mathbb R^m$ to the unit sphere in the Hilbert space, such that
$$curv^\perp_{ver_\infty}(\mathbb R^m\hookrightarrow S^\infty)=\sqrt \frac {(m-1)(2m+1)}{(m+1)^2}= 
\sqrt {2-\frac {5}{m+1}+\frac{2}{(m+1)^2}},$$}
where  equivariance is understood with respect to a certain unitary representation of the isometry group of $ \mathbb R^m$. 

{\color {red!50!black}Probably,} one can show that this $ver_\infty$ realizes the {\it  minimum} of the curvatures   among all equivariant  maps $\mathbb R^m\to S^\infty$.\vspace {1mm}

(b) Instead of   $ver_s$,  one could achieve (essentially) the same result with a use of compositions of the classical Veronese maps, $ver: S^{m_i}\to S^{m_{i+1}}$,  $_{i+1} = \frac {(m_i+1)(m_i+2)}{2}-2$,
$$S^{m_1} \hookrightarrow  S^{m_2}\hookrightarrow  ...  \hookrightarrow S^{m_i},$$
starting with   $m_1=2m-1$ and going up to $i=m$. (Actually, $i\sim \log m$ will do.)

\subsection { Product Manifolds, Connected Sums and  Related Constructions} 

Let  $f_i: X_i^{m_i}\to \mathbb B^{m_i+1}(1) $, $i=1,...,l$,  be immersions with focal radii $r$ 
 and let $f_0: X_0^{m_0}\to B^{l}(1) $ be an immersion with
 $foc.rad_{f}(X^{m_0}_0)= r_0$,

 Then the $\rtimes$-construction (see 4.1)  delivers  an immersion
$$f:X =\bigtimes_{0}^l X_i\to B^N(1), \mbox { }  N =l+\sum_1^l m_i,$$
such that
$$foc.rad_{f_\rtimes}(X_\times)\geq\max_{0<\lambda \leq1}  \frac{\min (r-\lambda, \lambda r_0)}{\sqrt l +\lambda r_0.} $$

Similarly,  if $X_0^{m_0}$ admits a  
$\nabla^\perp$-trivial (see 4.1) immersion to $B^M(1)$  
with focal radius $r_0$, then
$X$ admits an immersion to $B^{M+k}(1)$ for all $k \geq 1-M +
\sum_{0}^lm_i,$
such that
$$foc.rad_{f_\rtimes}(X_\times)\geq \max_{0<\lambda \leq1}  \frac{\min (r_0-\lambda, \lambda r/\sqrt l)}{\sqrt l + \lambda r/\sqrt l}$$

\textbf {6.2.A. Example: Product of Spheres.} Let  
$$X=X^m=\bigtimes_i S^m_i, \mbox { } \sum_im_i=m,$$
 and let $\mu=\min_im_i$. Then there exists an immersion $f:X\to B{^m+1}(1)$, such that 
$$curv_F^\perp(X)\leq const_\mu m^\frac {\mu+2}{\mu+1}$$

{\it Proof.} Adopt the torus-by-torus construction  \textbf {4.1.C} to product of spheres, where  instead  of  squaring   maps at each step, use (Cartesian) product of at least $\mu$ of maps, where  then  the above inequality for $foc.rad$ translated to curvature  apply.

{\it Embedding  Remark.} Observe that  the  resulting maps $X^m\to B^{m+1}(1)$ are {\it embeddings.}

\textbf {6.2.B. Connected Sums.} {\sf  If $m$-manifolds  $X_i$, $i=1,2,...,l$, admit  immersions to the unit ball   $B^n=B^n(1)$, $n>m$, with the  curvatures bounded by a constant $C$,  then the connected sum $X_1\#...\#X_l$  can immersed to  $B^n$  with curvature 
bounded by $5C$}.

{\it Proof}.  Make  {\it geometric connected  sums} of all  
 $X_i\hookrightarrow B^n$ with  the unit equatorial  sphere $S^m\subset S^n=\partial B^n$, where this is done  with 
each $X_i$   individually  with a copy of $S^m\subset B^n$ by connecting $X_i $  with $S_1^m=S^m$  with a tube with curvature $< 5C$.
Then the connected sum between $X_i$ is implemented by making similar tubes between 
$S^m_i$.

{\it Example.} Since there are 2-Tori in the unit 3-ball with $curv=3$, the minimal possible curvatures of orientable surfaces  $X$ satisfy
$$min.curv^\perp(X_{ori}^2\hookrightarrow B^3(1))<15,$$
while    non-orientable ones have 
$$min.curv^\perp(X^2\hookrightarrow B^3(1))\leq 5min.curv^\perp(
\mathbb RP^2\hookrightarrow B^3(1)) < 50,$$
the Boy surface seem to have  curvature about 10, 
Probably, all   surfaces have $min.curv^\perp<10$, but it is unclear, not even for the 2-torus, what actually 
 minimal curvatures of surfaces in $B^3(1)$ are.

\textbf {Attaching $k$-Handles for $k\geq 2$.} To  attach a  handle to 
a   sphere $ S^{k-1}\subset X$  with a controlled   the curvature, 
  with a controllable increase of  the  curvature,
one  needs   a regular  $\delta$-neighbourhood of this sphere  in $X$
  with $\delta$   controllably bounded from below: this which would allow attaching a $k$ handle   with the  curvature  increase roughly by  $1/\delta$.

 For instance, if $k=2$ an $S^1\subset X$  is the shortest non-contractible
curve in $X$, then it does admits such a  neighbourhood in $X$  with $\delta$   controllably bounded from below by the curvature of $X$; thus   attaching  with certain  normal frames 
2-handles  to it is possible with curvature increase by  a definite multiplicative constant.

In general one can show the following. 

\textbf {6.2.C.  Handles Stretch  Proposition.} (Compare with 4.3.C.)
{\sf Let an immersed manifold  $X_\diamond^m\overset {\phi}\hookrightarrow B^N(1)$ be obtained 
from $X^m \overset { f}\hookrightarrow B^n(1)$ by attaching $l$-handles for $l\leq k$ 
where,  all steps surgery  keep in the {\it class of {immersed manifolds}.

Then $\phi$ is    regularly homotopic to an immersion $\phi_1:X_\diamond \hookrightarrow B^n(1)$, such that
$$curv^\perp_{\phi_q}(X_\diamond) \leq C^{2k}curv_F^\perp(X)$$}
for $C\leq$ 10 000.}

{\it Sketch of the Proof.} Regularly  homotop  $f$ in $B^n(1)$ to an immersion $f_1$ 
with $curv^\perp_{f_1}(X)\leq 100^{2k} curv_F^\perp(X)$ and such that that the $f_1$-induced Riemannian metric  in a (small) neighbourhood $U$ of the  $2k$-skeleton of a  smooth triangulation of $X$ is by an arbitrarily large (independently of $U$)    factor $\lambda$ greater than the $f$-induced metric. 

Assume without loss of generality that all  spheres  $S^i$, at  which the surgery   performed  are located and in  $U$  don't  intersect there (this  is possible for $m\geq 2k$, which we may assume with no problem) and
choose $\lambda$ so large that the  union of   these spheres has a nice  thick regular neighbourhood,  where the surgery can be made with at most  
$100^{2k} $ increase in the curvature.

{\it Remark.} It is not hard to visualise      an actual proof along these lines but I don't  see how to write it down in a   readable form. 

\subsection {Embeddings with Small Curvatures}  

{\it Connected Sums of Embedded Manifolds.} {\it If $X=X^m$ admits an embedding (i.e. a immersion 
with no self-intersection)  to $B^{m+1}(1) $  with curvature $\leq c$, then 
 the connected sums of $2l$-copies of $X$ 
embed to $B^{m+1}(1)$ with curvatures $<100c$.}

{\it Proof.} Let  $X_1\subset B^{m+1}(1) $ be obtained from $X$ by attaching  a single 1-handle  
$S^{m-1} \times [0,1]$,  such that 
$curv X_1\subset B^{m+1}(1) < 10c $.

Let $ \tilde X_l$  be the natural cyclic covering of $X_1$ of order $l$ and let $\bar X_l$   be obtained
by cutting $ \tilde X_l$ along the sphere   $S^{m-1}\subset \tilde X_l$  from the handle.

Observe that this $\bar X_l$ is a manifold with two spherical boundary components  and that it 
 (almost) naturally embeds to $^{m+1}(1) $ with curvature $<10c$.

Let  $\bar X'_l\subset B^{m+1}(1)\setminus \bar X_l $ be obtained by a slight normal 
displacement of  $\bar X_l$ and let us attach  $\bar X'_l$ to  $\bar X_l$ along a pair of 
nearby $ (m-1)$- spheres and also fill in the remaining two boundary spheres with $m$-balls.
Clearly, the resulting manifold, call   it $X_{2l}$,  is diffeomorphic to the connected sum  
 of $2l$ copies of $X$  and it is not hard to arrange an  embedding of  $X_{2l}$ to the unit  ball 
 with curvature $<100$.

{\it Exercises.   } (a) Let $X=X^m$ be a connected sum of an arbitrary number  of  manifolds  
diffeomorphic to product of spheres.
 Show that $X$ embeds to the unit $(m+1)$-ball with curvature$<500\cdot2^\frac {m}{2}m^\frac {3}{2}$. 

 {\it Hint.} Embed mutually non-diffeomorphic  products of spheres into 
 $2^m$  disjoint   $r$-balls in $ B^{m+1}(1)$ of radii $r=2^{-\frac {m+2}{2}}.$


(b) Let $X=X^m$ be  
 disconnected closed manifold, which 
 contains $l$ {\it mutually non-diffeomorphic} components. Show that 
 $$curv_F^\perp(X\hookrightarrow B^{m+1})\geq const_m l,\mbox  { } const_n\geq \frac {1}{(10m)^m}, $$
 for all embeddings $f:X \hookrightarrow B^{m+1}(1)$.

(c) Construct closed m-dimensional manifolds $X_i$,   $i=1,  2,....$ for all $m\geq 6$,  such tat 
all of them embeds to 
$B^{7}(1)$   and such that  embedding   of  connected sums of
 $l$  among these  manifolds have curvatures $\geq const l$.

\textbf {Question.} Can one  have these  $X_i$  embeddable  to $\mathbb R^{m+1} $  with curvatures
< 1 000 000?


\subsection  {Cycles with  Small Curvature}

Our  equidimensional  expanding maps  are effective in delivering 
immersed submanifolds with controllably bounded curvatures, because these maps 
themselves, besides  being expanding, have controllably bounded second derivatives.

In general, it is {\color {red!60!black} hard} to \vspace {1mm}

{\sf  construct a immersed $m$-dimensional submanifolds $X \hookrightarrow Y$ with {\it small curvature}  and with {\it non-zero}  homology classes $[X]\in H_m(Y)$}.
 \vspace {1mm}

 Apparently, all known results  of this kind badly depend on the dimension and/or codimension
  of $X$, see [CDM2016]

A happy exception is the codimension one case,  $m=n-1$, where there is no topological obstructions for the existence of $X$ and where an   equidistant smoothing delivers hypersurfaces with controllably small curvatures as follows.\vspace{1mm}

 Let $Y$ be a {\it proper  Riemannian band} of dimension  $n$, that   is a Riemannian manifold, the boundary $\partial Y$ of which is divided into two disjoint parts, 
  $\partial Y=\partial_-Y\sqcup \partial_+ Y$, where $\partial_\pm Y$ are unions of connected components of $\partial Y$, and  denote by $d$ the {\it width} of $Y$,
$$d=width (Y) =_{def} dist(\partial_-Y, \partial_+ Y).$$

{\sf Let us $d_1$-equidistantly push $\partial_-Y$ inside $Y$ for $d_1< d$ and then  $d_2$-equidistantly move the resulting hypersurface, denoted  $\partial_{-d_1}$, back
 toward $\partial_-Y$   with $d_2<d_1$. }

That is,  $\partial_{-d_1}$ is equal to the (topological) boundary of the $d_1$-neighbourhood 
$U_{d_1}(\partial_-Y)\subset Y$ and the result of the second move, call it 
 $X_\circ=\partial_{-d_1|+d_2}\subset U_{d_1}(\partial_-Y)$, is the boundary of 
 $U_{d_2}(\partial_{-d_1})\subset U_{d_1}(\partial_-Y).$ 

Let us evaluate the curvature of $X_\circ$ in terms of the sectional curvatures of $Y$, where we  observe the following.

1. If $Y$   has constant sectional curvature $\pm\kappa^2$, then  $X_\circ$ is $C^{1,1}$-smooth and
$$foc.rad(X_\circ)\geq (\min (d_2, d_1-d-2));$$ 
accordingly  $curv^\perp(X)\leq \alpha^\pm_\kappa(\min (d_2, d_1-d-2))$
for the function $\alpha^\pm$ from 1.B.

2. If.more generally,  the sectional  curvatures  of  $Y$ is pinched between two values, that are the  curvatures of two standard surfaces $S_\pm$ with constant curvatures,
$$sect.curv^\perp (S_-)\leq sect.curv^\perp(Y)\leq sect.curv^\perp (S_+),$$ 
then the curvature of $X_\circ$ is bounded by the  maximum the two numbers:

$\bullet_1$ the first number is the  curvature of the circle of the radius $d_2$ in $S_-$;

$\bullet_2$ the second  number is the curvature of the circle $S^1(r) \subset S_+$,   such
that  the curvature of the  concentric circle $S^1(r+d_2)$  is equal  to the curvature of  the
 $d_1$-circle in $S_-$; 

It follows, for instance, that\vspace {1mm}

{\color {blue}({\Large $\circ$}$_d$)} {\sf if $$-1\leq sect.curv^\perp(Y)\leq 1$$
 and $d=width(Y)\leq 1$, then }
 
 {\it $Y$ contains a smooth hypersurface,  which {\it separates 
 $\partial_-Y$ from $\partial_+Y$ and such that}
 $$curv^\perp(X)\leq \frac {4}{d}.$$}
 \vspace{1mm}
 {\it Corollary.} {\sf  Let $Y$ be a complete Riemannian  $n$-manifold with 
 $|sectcurv^\perp(Y)|\leq \kappa^2$ and with $inj.rad(Y)\geq r$.}

 Then 
 
{\color {blue}({\Large $\circ$}$_{\kappa,r}$)} {\it all integer $(n-1)$-dimensional homology classes  $h\in H_{n-1}(Y)$ are realizable by smoothly immersed  oriented hypersurfaces $X\hookrightarrow Y$ with $curv^\perp(X)\leq 10\kappa+\frac {10}{r}$.}\footnote{If $Y$ is,  Riemannian flat, then the term $10/r$ is unneeded and if $Y$ is almost  flat one can do without  it  for  multiples of $h$  and     I am not certain about examples where the term $10/r$ is truly  needed. } 

\vspace {1mm}
Indeed, given a homology class $h\in H_1(Y)$, apply  {\color {blue}({\Large $\circ$}$_d$)}
 to the infinite cyclic covering of $Y$, which is defined by this  class. 

\vspace {1mm}

{\it Questions.} (a) Do {\color {blue}({\Large $\circ\circ$}$_d$)} and 
({\Large \color {blue}$\circ$}$_r$) meaningfully  generalize to
submanifolds $X\subset Y$ of codimensions $k>1$, where $Y$ is, in some way,  "wide  in $k$-directions"?

For instance, 
Let $Y$ be a Riemannian manifold homeomorphic  to $X_0\times B^k(1)$, where $X_0$ is a closed connected  orientable  manifold of dimension $n-k$, let the sectional curvature  of $Y$ be bounded 
by  $|\kappa(Y)|\leq 1$ and the injectivity radius by $inj.rad(Y)\geq 1$.\vspace {1mm}

{\sf What else need you  know about  $Y$ to effectively bound  the minimal possible  curvature  of a submanifold $X\subset Y$ {\it homologous to} 
$X_0=X_0\times \{0\}\subset X_0\times B^k(1)=X$?

What is the best bound on this curvature  in a presence of a {\it proper} (boundary-to-boundary)  {\it $\lambda$-Lipschitz} map $X\to B^k(1)$?\vspace {1mm}

Are, similarly to   {\color {blue}({\Large $\circ\circ$}$_\kappa,r$)},  non-zero multiples of the   homology classes $h\in H_m(Y)$, for all $m\leq dim(Y)$,
  realizable by immersed $m$-dimensional submanifolds 
$X\hookrightarrow Y$ with $curv^\perp(X)\leq 100m^{100}\left (\kappa+\frac {1}{r}\right)$? }
\vspace {1mm}

\textbf {From Focal Radius to Expansion.} Let us turn to the 

{\it opposite problem}:
{\sf In what cases does the   the $r$-neighbourhood $U_r(X)\subset X$ of an embedded manifold $X\subset Y$ with {\it "large" universal covering}, e. g.   for $X$ homeomorphic to $\mathbb T^m$,  
and   with large $ foc.rad(\mathbb X)$ 
  receive an  expanding map from a {"\it large manifold}" e.g.   from 
  $B^m(R)\times B^{n-m} \left (\frac {r}{100}\right)$  with large $R$?}

Here the answer is positive for $m=n-1$ and $m=n-2$: \vspace{ 1mm}

{\sl if $X$  receives expanding maps from the balls $B^m(R)$ for all $R$  (as e.g. the m-torus does), then, in the case  $m=n-1$, the neighbourhood  
 $U_r(X)$
receives expanding  maps from $B^m(R)\times B^{1}\left 
(\frac {1}{\sqrt 2} r-\varepsilon\right)$ for all  $R\to \infty $ and positive $ \varepsilon\to 0$.

And if $m=n-2$, then  $U_r(X)$ receives such maps from  $B^{m+1}(R)\times B^{1}\left (\frac {r}{2\sqrt 2}-\varepsilon\right)$.}
\vspace{ 1mm}

{\it Proof.}  The required map for $m=n-1$ and coorientable $X\subset Y$  is obtained with the obvious splitting 
$U_r(X)=X\times B^1(r)$ and the case $m=n-2$ follows by applying this to the  hypersurface 
$Z=\partial U_{r/2}(X)\subset U_r(X) $, where, clearly,  
$foc.rad(Z)=\frac {1}{2} foc.rad (X) \geq \frac {r}{2}$, and where the case of a  non-trivial 
normal bundle of $X\subset Y$ needs a little  thinking about.

\vspace {1mm}

But when it comes to $m\leq n-3$  nothing of the kind  seems to be true, where the apparent difficulty stems from the following phenomenon.

If $m, k\geq 2$, then the  topologically trivial sphere bundle $V=\mathbb R^m\times S^k\to \mathbb R^m$  admits an orthogonal connection $\nabla$ with an arbitrary small curvature 
such that all smooth  sections $\phi : \mathbb R^n\to V$ satisfy.
$$\sup_{x\in \mathbb R^m} ||\nabla \phi(x)||=\infty.$$

Despite this, our $U_r(X)$, still looks large for all $m$ and large $r=foc.rad(X)$,   but 
I don't know, how to make precise sense of largeness for these $U_r$.

Here is a specific question.

Let us regard $U=B^k(r)\times B^m(R)$ as (the total space of) a   $B^k(r)$-bundle  
over the ball $B^m(R)$,  let $\nabla$  be a Euclidean connection in this bundle and 
$g_\nabla$ the corresponding Riemannian metric on $U$, that is the sum of the differential quadratic form induced by the map $U=B^k(t)\times B^m(R)\to B^m(R)$ with 
 the  Euclidean metrics in the fibers $B_x^k(r)\subset U$, $x\in  B^m(R)$ extended to 
 $T(U)$   by zero on the $\nabla$-horizontal vectors.
\vspace {1mm}

{\sf For which  $r$, $ R$ and $\underline R$  the manifolds  
$ (U,g_\nabla)$ admit  {\it \color {red!40!black}no } expanding maps $ (U,g_\nabla)\to B^{m+k}(\underline R)$
for all connections $\nabla$?

Conversely, from what kind of manifolds do   $ (U,g_\nabla)$ receive expanding maps.?}

\section {References}

[BB2009] Eiichi Bannai,  Etsuko Bannai,
{\sl  A survey on spherical designs and algebraic combinatorics on spheres.}
European Journal of CombinatoricsVolume 30, 2009 pp 1392-1425. \vspace {2mm}

[Be2022]  S.M. Berge  {\sl Eigenvalues on Spherically Symmetric Manifolds,} 

arXiv:2203.11911.
    \vspace {2mm}

[Bet1991]  F. Bethuel, The approximation problem for Sobolev maps between two, Acta
Math. 167 (1991), 153–206.
\vspace {2mm}

[BHJ]A generalization of Geroch's conjecture
Simon Brendle, Sven Hirsch, Florian Johne {\sl A generalization of Geroch's conjecture}
	arXiv:2207.08617\vspace {2mm}



[CDM2016] Gregory R. Chambers, Dominic Dotterrer, Fedor Manin   and Shmuel Weinberger,
{\sl Quantitative null-cobordism},
Journal of the American Mathematical Society 31(4).

\vspace {2mm}

[CMS]        Otis Chodosh, Christos Mantoulidis, Felix Schulze.
{\sl Generic regularity for minimizing hypersurfaces in dimensions 9 and 10} arXiv:2302.02253v2.\vspace {2mm}

[DW1971]   M. P. do Carmo and  N. R. Wallach, {\sl Minimal immersions of spheres into spheres}. Ann. of Math. (2)93, 43–62 (1971).\vspace {2mm}

[FLM 1977] T. Figiel, J. Lindenstrauss, V. Milman, {\sl The dimension of almost spherical sections},
of convex bodies, Acta Math. 139 (1977), 53–94. \vspace {2mm}

[Ge 2021] J. Ge,  {\sl Gehring Link Problem, Focal Radius and Over-torical width},
arXiv:2102.05901 [math.DG] \vspace {2mm}

[GE1971] M. Gromov, Y. Eliashberg,  {\sl Removal of singularities of smooth maps}, Izv. Akad. Nauk, S.S.S.R.
35, \#5, pp. 600-627. \vspace {2mm}

[GE1971'] M. Gromov, Y. Eliashberg, {\sl Construction of nonsingular isoperimetric films}. Proceedings of the Steklov Institute of Mathematics, 1971, 116, 13-28. \vspace {2mm}

[Gr1986]   M. Gromov, {\sl  Partial differential relations}  Springer-Verlag (1986).  \vspace {2mm} 

[Gr1988] M. Gromov, {\sl  Dimension, non-linear spectra and width}, Lect. Notes in Math. Springer-Verlag 1317 (1988), 132-185. \vspace {2mm}

 [Gr2017]  M. Gromov,  {\sl Morse Spectra, Homology Measures and Parametric Packing Problems},
 arXiv:1710.03616

 \vspace {2mm}

 [Gr2018] M. Gromov,  {\sl Metric Inequalities with Scalar Curvature}
Geometric and Functional Analysis volume 28, pages 645–726 (2018) \vspace {2mm}

[Gr2021] M. Gromov,  {\sl  Four Lectures on Scalar Curvature},

arXiv:1908.10612.  \vspace {2mm}

   [Gr2022] M. Gromov,  {\sl Scalar Curvature, Injectivity Radius and
Immersions with Small Second Fundamental
Forms, } arXiv:2203.14013.  \vspace {2mm}

 [Gr2022'] M. Gromov,  {\sl  Isometric Immersions with Controlled Curvatures}, arXiv:2212.06122. 

 \vspace {2mm}

 [Gr2023]  M. Gromov,{\sl Product Inequalities for $\mathbb T^\rtimes$ -Stabilized Scalar Curvature}, arXiv:2306.02932.

  \vspace {2mm}

[Hit1974]  N. Hitchin, {\sl Harmonic Spinors}, Advances in Mathematics, 14(1), 1-55 (1974).\vspace {2mm} 

[K1995] H. Konig, {\sl Isometric imbeddings of Euclidean spaces into finite  dimensional $l_p$ -spaces, }
Banach Center Publications (1995)
 Volume: 34, Issue: 1, page 79-87. \vspace {2mm}

 [Ku2021]  Y.Kubota,   {\sl Band width and the Rosenberg index}
arXiv:2108.08506 [math.KT]
  	 \vspace {2mm}

[LB2021] C. LeBrun, {\sl On the Scalar Curvature of 4-Manifolds}

arXiv:2105.10785 [math.DG] \vspace {2mm}

LW1993] Yu. Lyubich, L. Vaserstein, {\sl  Isometric imbeddings between classical Banach spaces,
cubature formulas, and spherical designs,} Geom. Dedicata 47 (1993), 327–362. \vspace {2mm}

[PVZ 2017] G. Paouris, P. Valettas and J. Zinn,   {\sl   Random version of Dvoretzky’s theorem in $l_p^n$, }  	arXiv:1510.07284 


 \vspace {2mm} 

[Pet 2023]
A. Petrunin {\sl Gromov’s torii are optimal} arXiv:2304.00886v1 \vspace {2mm} 

 [Sm1993]   N. Smale,
{\sl Generic regularity of homologically area minimizing hyper
surfaces in eight-dimensional mani-
folds}, Comm. Anal. Geom. 1, no. 2 (1993), 217-228.\vspace {1mm}
  \vspace{1mm}

[Stolz(survey) 2001] S. Stolz {\sl Manifolds of Positive Scalar Curvature}  \url {http:
//users.ictp.it/~pub_off/lectures/lns009/Stolz/Stolz.pdf} \vspace {2mm}

[SY1979] R. Schoen and S. T. Yau, {\sl On the structure of manifolds
with positive scalar curvature}, Manuscripta Math. 28 (1979).\vspace {2mm} 
159-183.

[SY2017] R. Schoen and S. T. Yau, {Positive Scalar Curvature
and Minimal Hypersurface Singularities}, arXiv:1704.05490.\vspace {2mm}

[WXY2021] J. Wang, Z. Xie, G. Yu,   {\sl  An index theoretic proof of Gromov’s
cube inequality on scalar curvature}, ¯arXiv:2105.12054

\end{document}